\documentclass[10pt]{amsart}

\usepackage{amsmath,amsthm,amssymb,epsf}
\usepackage[hidelinks]{hyperref}
\usepackage{graphicx,epstopdf}
\usepackage[margin=1in]{geometry}

\newtheorem{theorem}{Theorem}[section]
\newtheorem{proposition}[theorem]{Proposition}

\theoremstyle{definition} \newtheorem{definition}[theorem]{Definition}
\newtheorem{example}[theorem]{Example}
\newtheorem{construction}[theorem]{Construction}

\theoremstyle{remark} \newtheorem{remark}[theorem]{Remark}

\DeclareMathOperator{\fr}{frame}
\DeclareMathOperator{\KS}{KS}
\DeclareMathOperator{\planar}{planar}

\newcommand{\Hc}{\mathcal{H}}
\newcommand{\R}{\mathbb{R}}

\title[Two types of Heegaard diagram]{On two types of Heegaard diagram used in knot Floer homology}
\author{Andrew Manion}

\begin{document}

\begin{abstract} In this brief note, we give an explicit sequence of Heegaard moves interpolating between local versions of the Kauffman-states Heegaard diagram and the planar Heegaard diagram used in knot Floer homology, and show how these local moves can be used to go between the global versions of the Heegaard diagrams.
\end{abstract}

\maketitle

\section{Introduction}

Knot Floer homology \cite{HFKOrig, RasmussenThesis} is a powerful modern invariant for knots and links in $3$-manifolds. It is related to the Heegaard Floer invariants for $3$- and $4$-manifolds defined by Ozsv{\'a}th--Szab{\'o} \cite{HFOrig, HF4Mfd}, it categorifies the Alexander polynomial, and it exactly detects many important topological properties of knots that are only bounded or restricted by the Alexander polynomial.

To define the knot Floer homology of a knot or link $L$, one must first choose a Heegaard diagram representing $L$; in general this can be a multi-pointed diagram as in \cite[Section 2]{MOS}. From the diagram, one gets a graded free abelian group with a natural basis; after making additional analytic choies, one gets a differential on this group such that the homology of the resulting complex (up to isomorphism, and after correcting for tensor factors that result from having multiple basepoints) is an invariant of $L$. 

There are always infinitely many diagrams representing any knot or link. In order to study the algebraic properties of knot Floer homology, various systematic procedures for constructing multi-pointed Heegaard diagrams from knot projections have been considered. 

One such procedure, introduced in \cite{AltKnots}, constructs a two-pointed Heegaard diagram from a projection of a knot or link in $S^3$. Generators of the chain complex associated to this diagram are in bijection with Kauffman states \cite{FKT} of the projection. This family of Heegaard diagrams, the ``Kauffman-states diagrams,'' have recently provided \cite{OSzNew, OSzNewer, OSzHolo, OSzPong} the basis for a deep and sophisticated algebraic apparatus for knot Floer homology in the spirit of bordered Floer homology \cite{LOTBorderedOrig}, raising hopes for new relationships between knot Floer homology and constructions in representation theory (see \cite{ManionDecat, ManionKS}). This procedure was generalized by Ozsv{\'a}th--Stipsicz--Szab{\'o} in \cite{OSSz} to define diagrams for singular as well as nonsingular crossings.

Another procedure for constructing multi-pointed Heegaard diagrams from knot projections yields the ``planar'' Heegaard diagrams. These diagrams were first defined for fully singular link diagrams in \cite[Section 5]{OSSz} and generalized for all link diagrams in \cite[Figure 12]{OSzCube} (the exact version we consider was defined in \cite{ManolescuCube}). The planar diagrams seem especially suited for studying the relationship between Heegaard Floer homology and HOMFLY-PT homology; see \cite{DowlinComposition, DowlinKR, DowlinCategorification}.

Given the recent results about knot Floer homology based on both types of diagrams, it is natural to ask whether the two procedures mentioned above are related. In principle, any two multi-pointed Heegaard diagrams representing the same link can be related by a sequence of Heegaard moves, i.e. isotopies, stabilizations (of two types), destabilizations, and handleslides. The proof of this fact relies on Morse theory and is non-constructive; it is desirable to have a constructive proof that the specific diagrams in question are related by Heegaard moves.

Indeed, this question was addressed for fully singular link projections in \cite[proof of Lemma 3.7]{OSzCube}. Ozsv{\'a}th--Szab{\'o} construct a sequence of Heegaard moves transforming the Kauffman-states diagram into the planar diagram for a fully singular link projection. Besides the limitation to singular projections, however, another potentially undesirable feature of this construction is that it treats crossings differently depending on the position of their edges in an ordering on all the edges of the link projection. 

In this note, given a link projection that may have nonsingular as well as singular crossings, we construct an explicit sequence of Heegaard moves which transform the Kauffman-states diagram of the projection into the planar diagram of the projection. Our construction is more local in character than the one given in \cite[proof of Lemma 3.7]{OSzCube} for fully singular projections; all crossings (whether singular or nonsingular) are treated similarly and independently of the rest of the diagram. Our main result may be summarized as follows:

\begin{theorem}
The sequence of moves outlined in Construction~\ref{constr:Main} below transforms the Kauffman-states diagram of a (possibly singular) link projection into the planar diagram of the projection.
\end{theorem}

We use a version of the Kauffman-states diagram with ``ladybugs'' added as in \cite{ManolescuUnoriented}; see Section~\ref{sec:KSDiag} below. Ladybugs are often useful when working with the Kauffman-states diagram; see for example \cite{BaldwinLevine}. The diagrams with ladybugs can be transformed into the original diagrams using Heegaard moves, including index zero/three stabilizations, discussed in Section~\ref{sec:Moves} below. We review one way of doing this in Example~\ref{ex:LadybugMoves}.

Construction~\ref{constr:Main} will start by modifying a stabilized version of the Kauffman states diagram near each crossing, before performing additional moves that do not depend on whether the crossings are positive, negative, or singularized. To deal with the local constructions, we will introduce some basic definitions for Heegaard diagrams with boundary in Section~\ref{sec:Prelims}. The local constructions will be given in Section~\ref{sec:Local}; they will be patched together into moves relating the global Heegaard diagrams in Section~\ref{sec:Global}.

\subsection{Motivation}

One motivation for this work is to understand the planar diagram, and eventually its associated Heegaard Floer invariants, locally using the methods of bordered Floer homology. As a prerequisite for this algebraic setup, one needs to choose a decomposition of the Heegaard diagram into pieces satisfying certain properties. Combinatorial details of the pieces determine how the algebraic invariants work and how easy it is to define them. There are various possible decompositions of the planar diagram to which one could apply the ideas of bordered Floer homology; in this paper we study a decomposition related to one which has already been used for bordered Floer homology in \cite{OSzNew,OSzNewer}.

In particular, applying the bordered-Floer methods used in \cite{OSzNew,OSzNewer}, one should be able to define a bimodule over the algebras from \cite{OSzNew,OSzNewer} associated to each ``row'' of a Heegaard diagram like the one shown in the top-right corner of Figure~\ref{fig:MainConst}. One would then expect homotopy equivalences between these bimodules and the bimodules defined for nonsingular crossings in \cite{OSzNew} or \cite{OSzNewer} (depending on the chosen algebraic setup) and for singular crossings in \cite{ManionSingular} (in preparation); however, as currently planned, only local bimodules for a singular crossing between two strands will be defined in \cite{ManionSingular}.

One reason why the hypothetical bimodules of the previous paragraph would be interesting is that they could be used to give a local formulation of the singular skein exact triangle and cube of resolutions in knot Floer homology \cite{OSzCube, ManolescuCube}. The construction of this exact triangle uses the planar diagrams for a nonsingular crossing, its singularization, and its resolution. The corresponding bimodules should satisfy a mapping-cone relation. Constructing these bimodules and proving such a relationship could lead to local versions of results of Dowlin \cite{DowlinComposition, DowlinKR, DowlinCategorification} and further progress in the relationship between knot Floer homology and HOMFLY-PT homology.

\subsection{Acknowledgments}

I would like to thank Zolt{\'a}n Szab{\'o} for introducing me to the Heegaard diagrams behind \cite{OSzNew,OSzNewer}, in particular the stabilizations which proved essential in the constructions of this paper. I would also like to thank Nathan Dowlin, Aaron Lauda, Robert Lipshitz, Ciprian Manolescu, and Ina Petkova for useful conversations. This work was supported by the National Science Foundation under Grant No. DMS-1502686.

\section{Preliminaries}\label{sec:Prelims}

\subsection{Partial Heegaard diagrams}\label{sec:PartialHD}

Since we do not discuss algebraic invariants in this paper, the following (overly permissive in general) definition will suffice for our local constructions.

\begin{definition}
A partial Heegaard diagram is a compact oriented surface $\Hc$ with boundary, equipped with two finite collections of circles (the $\alpha$ and $\beta$ circles) smoothly embedded in $\Hc \setminus \partial \Hc$, two finite collections of compact intervals (the $\alpha$ and $\beta$ arcs) smoothly embedded in $\Hc$, and a collection of basepoints labelled either $X$ or $O$ which does not intersect $\partial H$ or the $\alpha$ or $\beta$ circles or arcs, such that:
\begin{itemize}
\item The endpoints of the $\alpha$ and $\beta$ arcs lie in $\partial H$, and these arcs intersect $\partial H$ transversely
\item There are no intersections between an $\alpha$ arc or circle and another $\alpha$ arc or circle, and similarly for $\beta$ arcs and circles
\item Any intersections between $\alpha$ and $\beta$ arcs or circles are transverse and occur in $H \setminus \partial H$. 
\end{itemize}
We will always consider partial Heegaard diagrams up to boundary-preserving isotopies that do not change any of the combinatorial properties of the basepoints or $\alpha$ and $\beta$ arcs and circles. As we will see below, we are even more interested in a coarser equivalence relation on partial Heegaard diagrams that (in particular) allows more general types of isotopy.
\end{definition}

We briefly review how Heegaard diagrams are used to represent $3$-manifolds and links. Suppose that $\partial \Hc = \varnothing$ and that adding $3$-dimensional $2$-handles to $\Hc \times [0,1]$ along the circles $\alpha_i \times \{1\}$ and $\beta_i \times \{0\}$ yields a $3$-manifold with boundary $S^2 \sqcup -S^2$ (in particular, the number of $\alpha$ circles and the number of $\beta$ circles must each equal the genus of $\Hc$). In this case, we may form a closed $3$-manifold $Y$ by filling in the $S^2$ boundary components with copies of $D^3$. The orientation on $\Hc$ allows us to specify an orientation on $Y$. The Heegaard diagram $\Hc$ (disregarding the basepoint data) is said to represent the closed oriented $3$-manifold $Y$. The Kauffman states and planar diagrams considered in this paper represent $S^3$.

For a Heegaard diagram $\Hc$ representing a closed $3$-manifold $Y$, assume that each connected component of the complement of the $\alpha$ circles in $\Hc$ contains exactly one $X$ basepoint and one $O$ basepoint, and that the same is true for the complement of the $\beta$ circles. If we first remove a small disk around each basepoint of $\Hc$, thicken $\Hc$ and attach handles as above, and omit adding copies of $D^3$, the resulting $3$-manifold is the complement of a link $L$ in the closed $3$-manifold $Y$ represented by $\Hc$. From the types ($X$ and $O$) of the basepoints, we can give $L$ an orientation (see e.g. Figure~\ref{fig:OrigPartialKSDiags}). We say that $\Hc$ represents $L \subset Y$. This definition can be generalized for singular links. The Kauffman states and planar diagrams associated to a projection of an oriented (and possibly singular) link $L \subset S^3$ represent $L$.

\begin{remark}\label{rem:Sutured}
More generally, if $\Hc$ has boundary but no $\alpha$ or $\beta$ arcs, one can form a sutured $3$-manifold from $\Hc$ by performing a similar thickening-and-handle-addition procedure. When $\Hc$ represents $L \subset Y$ as above and we construct a Heegaard diagram $\Hc^{\circ}$ from $\Hc$ by cutting out a small disk around every basepoint, then $\Hc^{\circ}$ has no $\alpha$ or $\beta$ arcs, so it represents a sutured manifold. This manifold is the link complement in which each component of the boundary has been given a meridional suture for each $X$ and $O$ basepoint corresponding to the component.
\end{remark}

\subsection{Gluing}\label{sec:Gluing}

Suppose $\Hc$ and $\Hc'$ are partial Heegaard diagrams. Let $\{C_1, \ldots, C_k\}$ and $\{C'_1,\ldots, C'_k\}$ be sets of connected components of $\partial H$ and $\partial H'$ respectively, with $k$ elements each, and suppose that for each $C_i$ we are given a diffeomorphism from $C_i$ to $-C'_i$ restricting to a bijection between the sets of endpoints of $\alpha$ and $\beta$ arcs intersecting $C_i$ and $C'_i$. We can form a new partial Heegaard diagram $\Hc \cup \Hc'$ by gluing each component $C_i$ to $C'_i$ via the given diffeomorphisms; we are omitting some details here, but these are not worrisome because we are only dealing with circles and surfaces.

\subsection{The Kauffman-states diagram}\label{sec:KSDiag}

Let $L$ be a (generic) projection of a link in $S^3$ onto the plane $\R^2$; for simplicity, assume that $L$ has a unique global minimum point in this plane. Also choose a finite set of points $p_i$ on $L$, away from the crossings and the global minimum. When we are given this data, we will call $L$ a marked projection.

Assume that every component of $L$ contains either the global minimum or at least one $p_i$. We will define the Kauffman-states Heegaard diagram for $L$ with the distinguished point placed at the global minimum and with ladybugs added at the points $p_i$. The choice of the $p_i$ (given the above restriction) is irrelevant up to Heegaard moves, which will be discussed in Section~\ref{sec:Moves} below (see Example~\ref{ex:LadybugMoves} in particular). 

\begin{figure}
\includegraphics[scale=0.75]{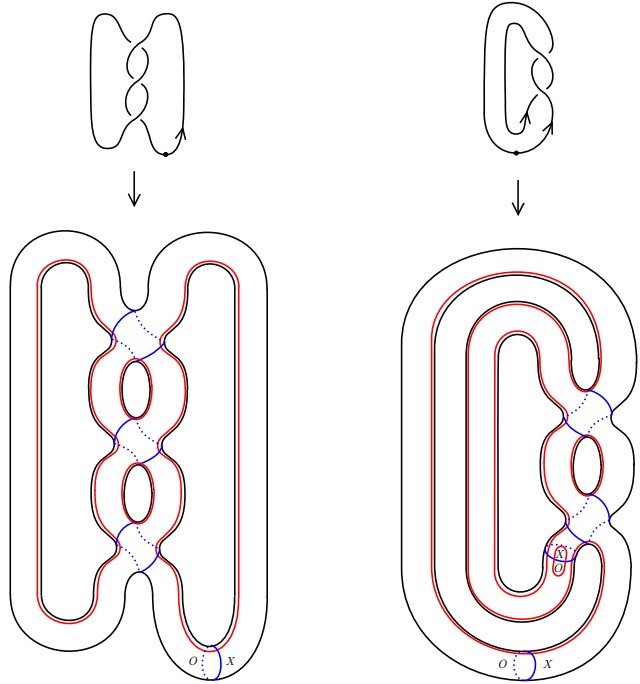}
\caption{Kauffman states diagrams for the trefoil and Hopf link.}
\label{fig:OrigKSDiagExamples}
\end{figure}

\begin{example} The Kauffman-states Heegaard diagrams for marked projections of the trefoil and Hopf link are shown in Figure~\ref{fig:OrigKSDiagExamples}. The global minimum of each projection is shown with a dot; the marked points $p_i$ are shown with dashes. 
\end{example}

\begin{figure}
\includegraphics[scale=0.75]{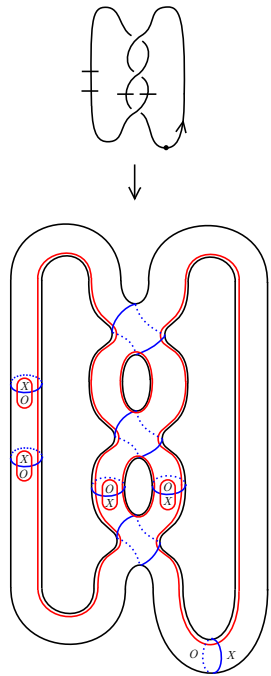}
\caption{Kauffman states diagrams for a projection of the trefoil with more markings.}
\label{fig:TrefoilMoreLadybugs}
\end{figure}

\begin{example} 
The Kauffman-states Heegaard diagram for a different marked projection of the trefoil, with more marked points (and thus more ladybugs), is shown in Figure~\ref{fig:TrefoilMoreLadybugs}.
\end{example}

\begin{figure}
\includegraphics[scale=0.75]{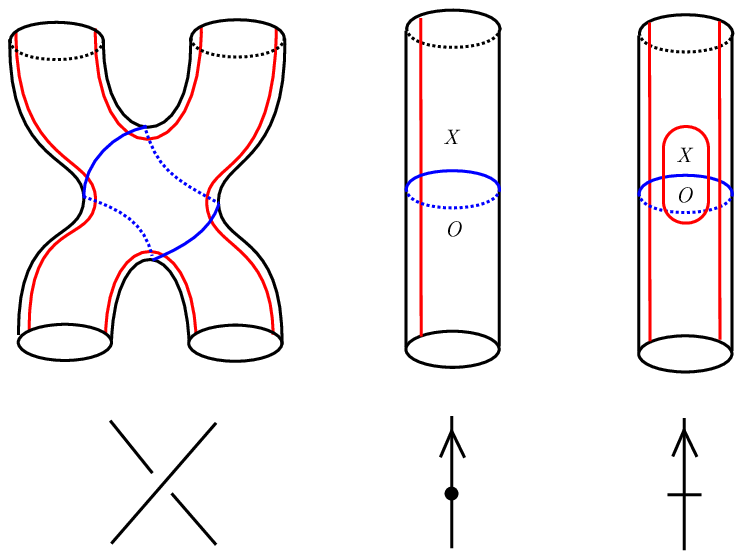}
\caption{The local pieces of the Kauffman-states diagram.}
\label{fig:OrigPartialKSDiags}
\end{figure}

In general, we can define the Kauffman states Heegaard diagram by gluing it together from the partial Heegaard diagrams shown in Figure~\ref{fig:OrigPartialKSDiags}. When an $\alpha$ arc in one of the local pieces would be adjacent to the unbounded region of the projection in the glued-together diagram, it should be omitted from the local piece. 

\begin{figure}
\includegraphics[scale=0.75]{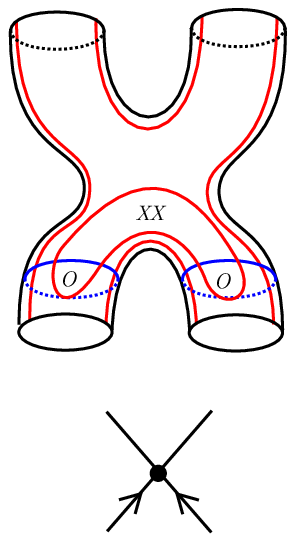}
\caption{The local piece of the Kauffman-states diagram for a singular crossing.}
\label{fig:OrigPartialSingKS}
\end{figure}

Following Ozsv{\'a}th--Stipsicz--Szab{\'o} \cite{OSSz}, we can also define a partial Heegaard diagram for a singular crossing, shown in Figure~\ref{fig:OrigPartialSingKS}. Diagrams for singular crossings are glued together into a global diagram just as for nonsingular crossings.

\subsection{The planar diagram}\label{sec:Planar}
We will also define the planar diagram from \cite[Section 4]{OSzCube} of a link projection $L$. This diagram is called planar because it has $\Hc = S^2$, unlike the Kauffman-states diagram. 

\begin{figure}
\includegraphics[scale=0.75]{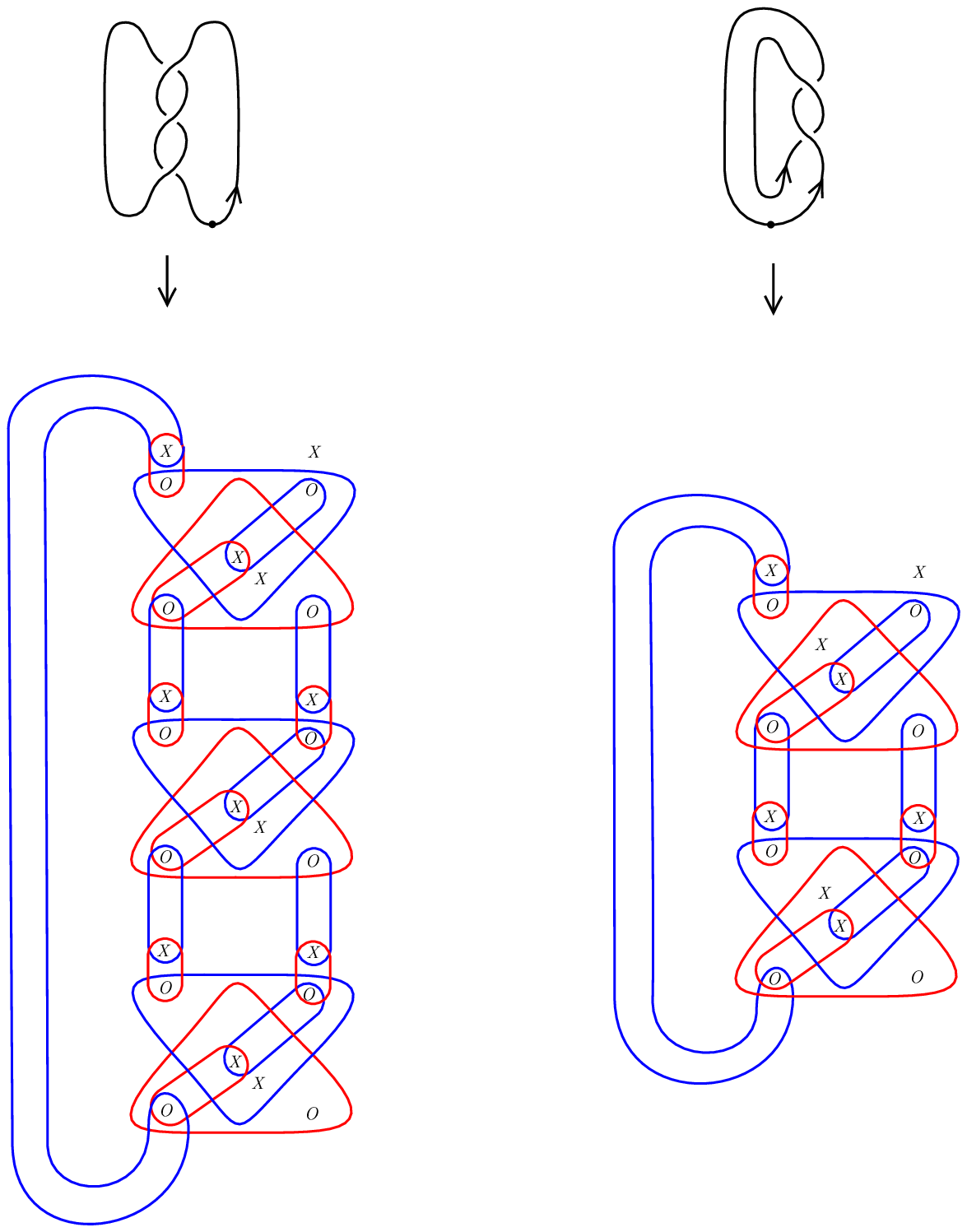}
\caption{Planar diagrams for the trefoil and Hopf link.}
\label{fig:OrigPlanarDiagExamples}
\end{figure}
\begin{example}
The planar Heegaard diagrams for projections of the trefoil and Hopf link are shown in Figure~\ref{fig:OrigPlanarDiagExamples}. Each picture should be interpreted as being embedded in a $2$-sphere which we will not draw.
\end{example}

\begin{figure}
\includegraphics[scale=0.75]{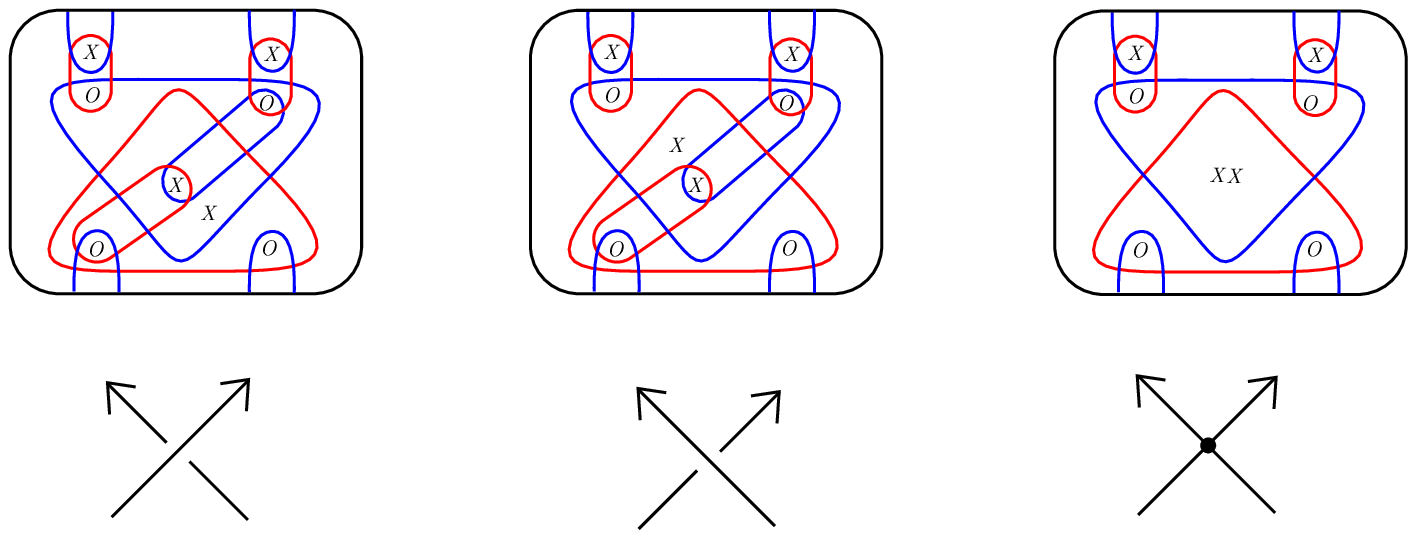}
\caption{Partial planar diagrams for each type of crossing.}
\label{fig:OrigPartialPlanarDiags}
\end{figure}

\begin{figure}
\includegraphics[scale=0.75]{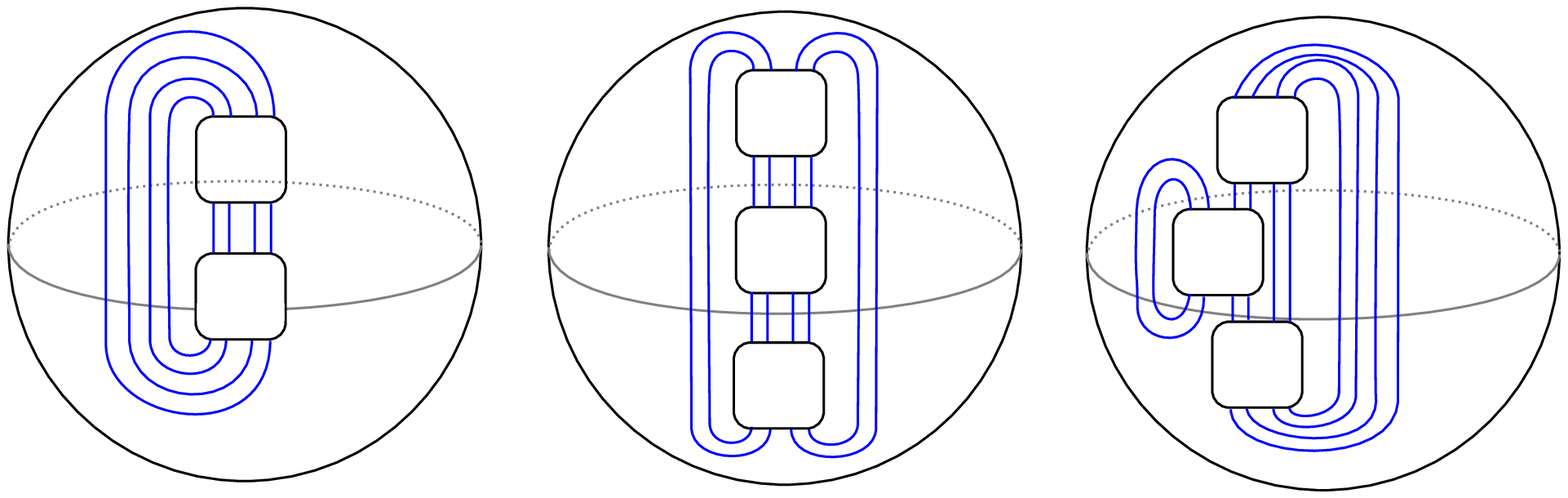}
\caption{Diagrams for forming the planar diagram from the diagram for each crossing: some examples.}
\label{fig:PlanarOutsidePieces}
\end{figure}

As with the Kauffman-states diagram, we can define the planar diagram by gluing together partial Heegaard diagrams. The partial diagrams for singular and nonsingular crossings are shown in Figure~\ref{fig:OrigPartialPlanarDiags}. To form the full planar diagram, one first glues the diagram for each crossing into an ``outside'' diagram depending on the projection, as in Figure~\ref{fig:PlanarOutsidePieces}. Then one removes the $\beta$ circle corresponding to the distinguished edge (the global minimum) and its corresponding small $\alpha$ circle.

\begin{figure}
\includegraphics[scale=0.75]{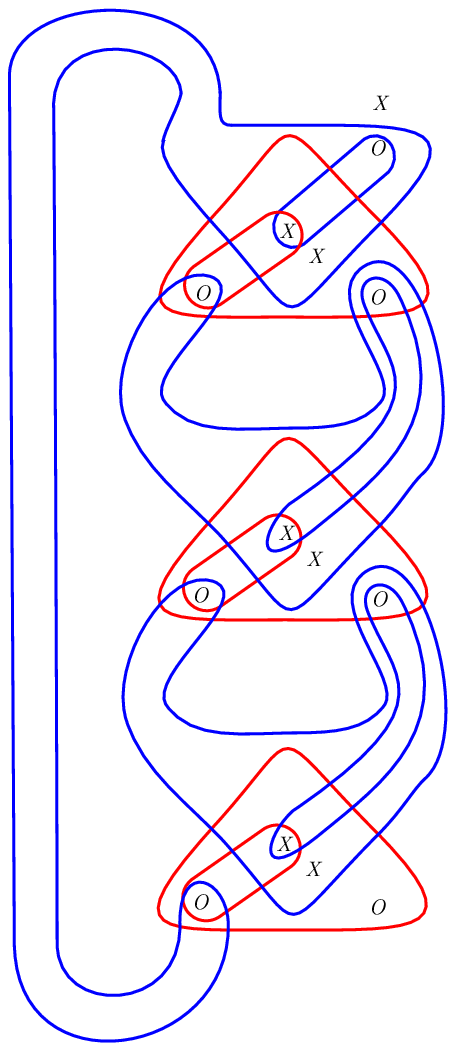}
\caption{A smaller variant of the planar Heegaard diagram for the trefoil.}
\label{fig:SmallerPlanarTrefoil}
\end{figure}

\begin{remark}\label{rem:SmallerPlanarDiags}
Often one works with a smaller variant of the planar diagram, shown for the trefoil in Figure~\ref{fig:SmallerPlanarTrefoil}. This smaller diagram can be obtained from the one we use by a series of handleslides and index zero/three stabilizations.
\end{remark}

\subsection{Heegaard moves}\label{sec:Moves}

We now discuss certain modifications one can perform to a partial Heegaard diagram, called Heegaard moves.

\begin{remark}
If a Heegaard diagram represents a closed $3$-manifold $Y$ or a link $L \subset Y$, this fact is invariant under performing Heegaard moves on the diagram. This is not quite true in the sutured setting of Remark~\ref{rem:Sutured}; the index zero/three stabilization described below does not leave the sutured $3$-manifold unchanged, but it modifies the sutured structure of the manifold in a predictable way.
\end{remark}

\begin{figure}
\includegraphics[scale=0.75]{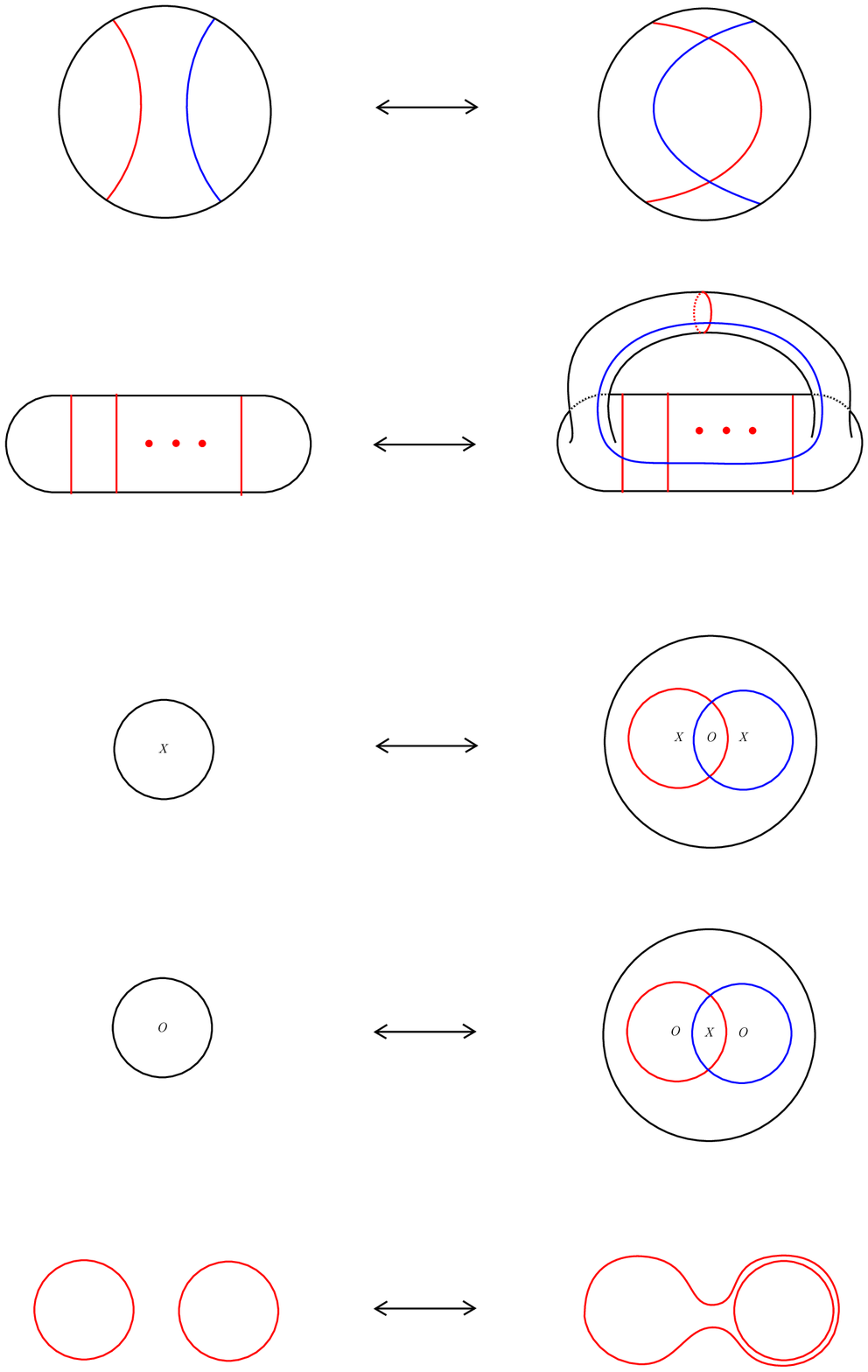}
\caption{Heegaard moves; all moves with $\alpha$ and $\beta$ interchanged are also valid.}
\label{fig:Moves}
\end{figure}

\begin{definition}\label{def:Moves}
For us, a Heegaard move is a modification of a partial Heegaard diagram in any of the following ways. A variation of this set of moves can be found in \cite[Proposition 3.3]{HFLOrig}; for the specific form of the index zero/three stabilization given below, see \cite[Lemma 2.4 and Figure 5]{MOS}. The terminology for stabilizations comes from Morse theory.
\begin{itemize}

\item First, one is allowed to isotope the Heegaard diagram as in the top row of Figure~\ref{fig:Moves}. These isotopies are allowed to change the combinatorial properties of the $\alpha$ and $\beta$ curves (e.g. their number of intersections), but the boundary of the Heegaard diagram must remain fixed and curves are not allowed to move past basepoints.

\item Second, one is allowed to perform an index one/two stabilization as in the second row of Figure~\ref{fig:Moves}, or to undo such a stabilization. The analogous stabilization with the roles of the $\alpha$ and $\beta$ curves reversed is also allowed. Note that the given form of the index one/two stabilization can be obtained from the simpler form described e.g. in \cite[Proposition 3.3]{HFLOrig} by using handleslides; we find it convenient to work with the given form.

\item Third, one is allowed to perform an index zero/three stabilization as in the third and fourth rows of Figure~\ref{fig:Moves}, or to undo such a stabilization. 

\item Finally, one is allowed to handleslide an $\alpha$ circle over an $\alpha$ circle, or a $\beta$ circle over a $\beta$ circle, as in the fifth row of Figure~\ref{fig:Moves}.
\end{itemize}
\end{definition}

\begin{figure}
\includegraphics[scale=0.75]{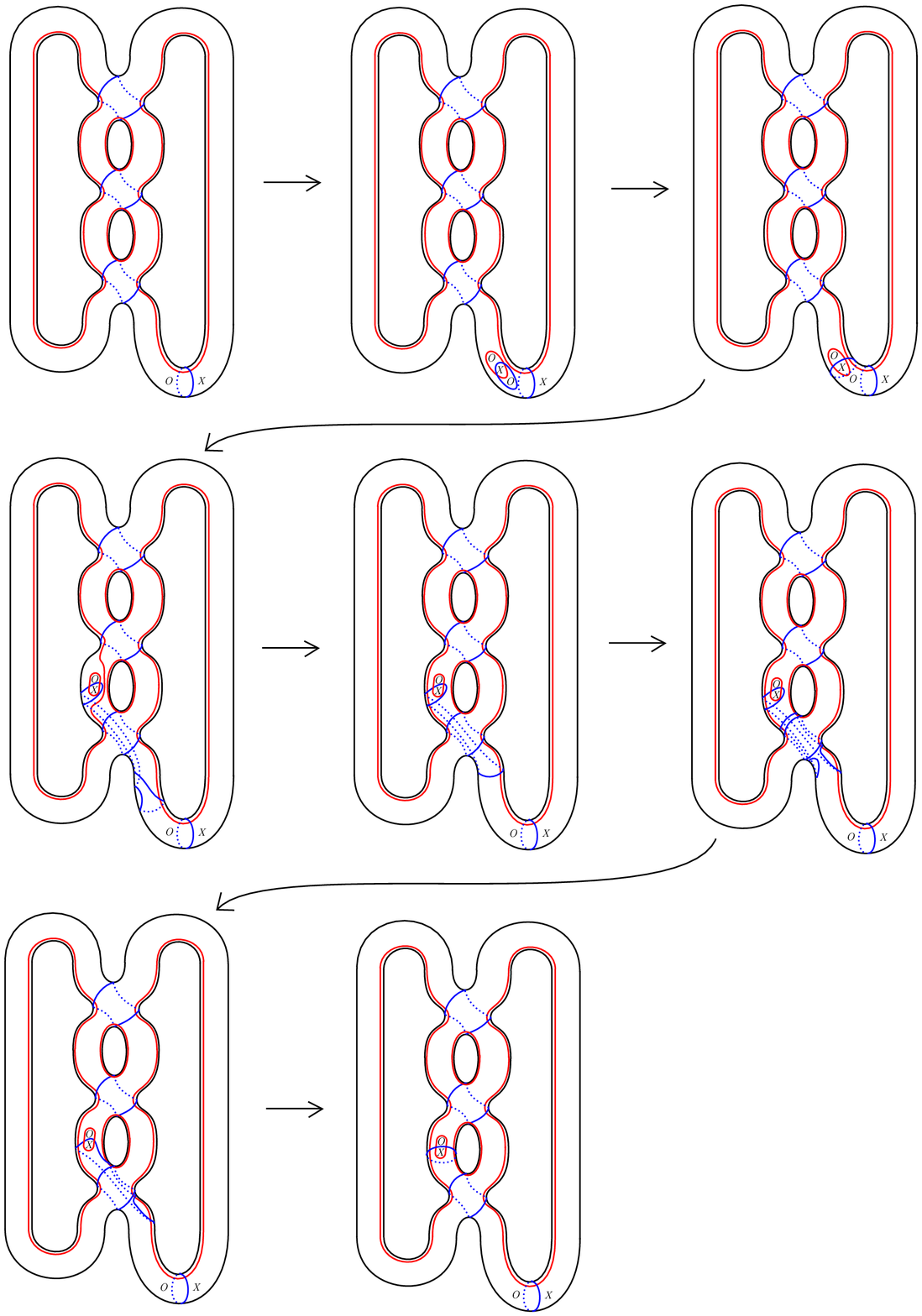}
\caption{Addition of a ladybug using Heegaard moves.}
\label{fig:Ladybug}
\end{figure}

\begin{figure}
\includegraphics[scale=0.75]{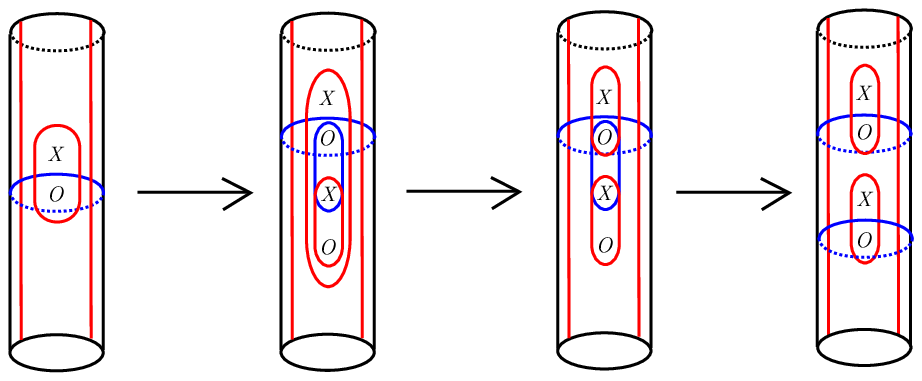}
\caption{Doubling a ladybug using Heegaard moves.}
\label{fig:LadybugDoubling}
\end{figure}

\begin{example}\label{ex:LadybugMoves}
In the case of the Kauffman-states diagram for the trefoil from Figure~\ref{fig:OrigKSDiagExamples} above, we show how to implement the addition of a ladybug using Heegaard moves in Figure~\ref{fig:Ladybug}. The first step is an index zero/three stabilization. The second step is a handleslide of a $\beta$ circle; at this point we have produced a ladybug near the global minimum of the projection.

The remaining steps show how to move a ladybug past a (nonsingular) crossing. The third step is a combinatorics-preserving isotopy to make things easier to see. The fourth step is a handleslide of an $\alpha$ circle, together with another combinatorics-preserving isotopy for simplification. The fifth step is a handleslide of a $\beta$ circle. The sixth and seventh step are isotopies that reduce the number of intersections between $\alpha$ and $\beta$ circles.

Ladybug addition works like this in general if $L$ is nonsingular. To add a ladybug on a component of $L$ that does not contain the global minimum, we start with an index zero/three stabilization at a pre-existing ladybug on this component, as shown in Figure~\ref{fig:LadybugDoubling}. 

If $L$ has singular crossings, one can create ladybugs at either of the two incoming edges of a singular crossing, just as one does at the global minimum. These ladybugs can then be propagated backwards along $L$, away from the singular crossing, until another singular crossing is reached which can create its own ladybugs. In this way, one can go from the Kauffman-states diagram without ladybugs to the Kauffman-states diagram for an arbitrary marking of the projection $L$ using the Heegaard moves of Definition~\ref{def:Moves}.

Adding or removing ladybugs is the only construction in this paper requiring the link projection to be traversed in-line as in \cite[proof of Lemma 3.7]{OSzCube}.
\end{example}

\subsection{Gluing and Heegaard moves}

The following fact will be useful in assembling our sequence of Heegaard moves from local pieces; in fact, our main justification for describing our construction as ``local'' is that this proposition is used repeatedly in Construction~\ref{constr:Main}.

\begin{proposition}\label{prop:Gluing} If $\tilde{\Hc}$ and $\tilde{\Hc}'$ are partial Heegaard diagrams obtained from $\Hc$ and $\Hc'$ by Heegaard moves, and we glue $\Hc$ and $\Hc'$ to get $\Hc''$ as in Section~\ref{sec:Gluing}, then we may also glue $\tilde{\Hc}$ and $\tilde{\Hc}'$ to obtain a partial Heegaard diagram $\tilde{H}''$, and the sequences of Heegaard moves transforming $\Hc$ and $\Hc'$ into $\tilde{H}$ and $\tilde{H'}$ may be combined to transform $\Hc''$ into $\tilde{H}''$.
\end{proposition}

\begin{figure}
\includegraphics[scale=0.75]{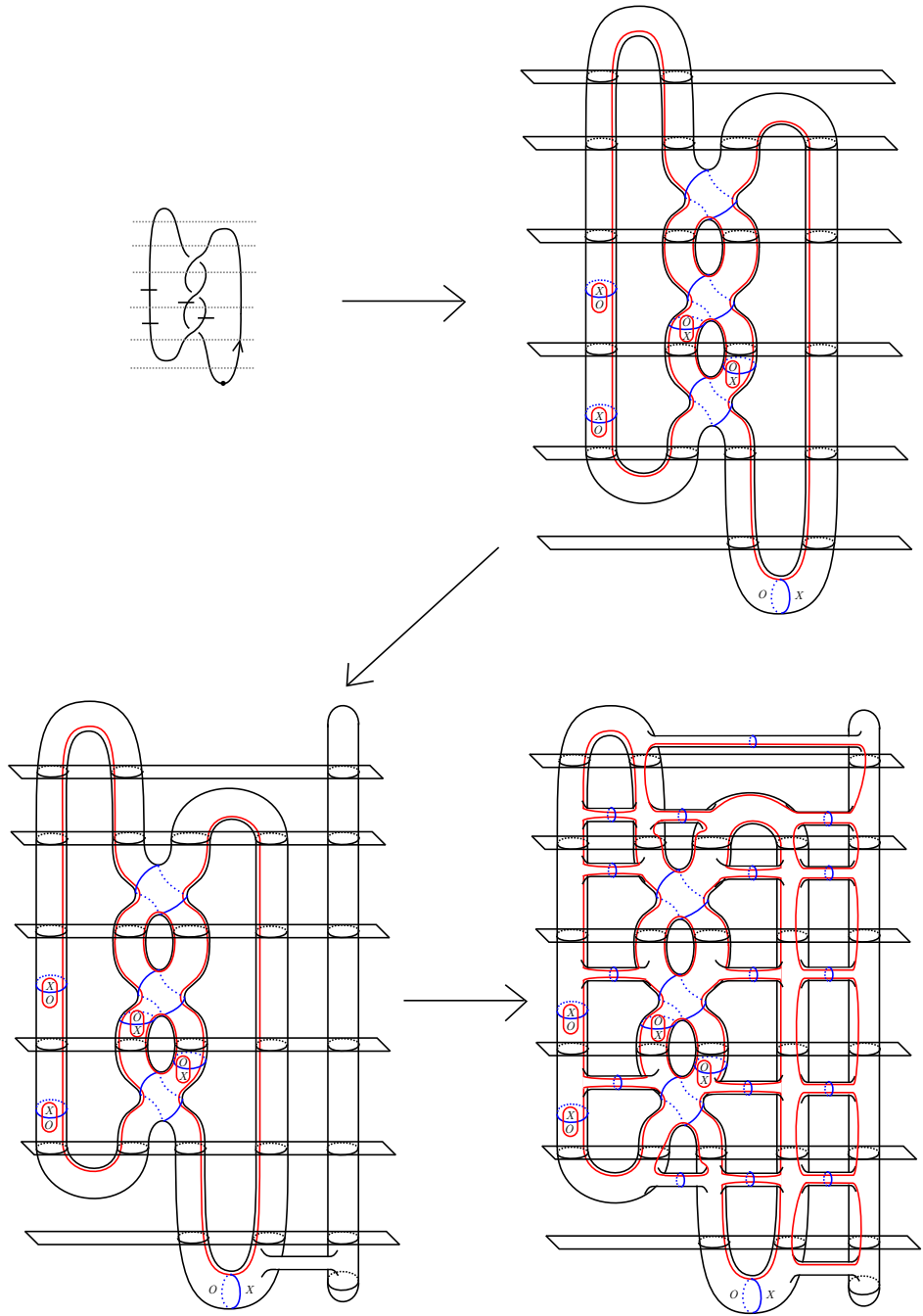}
\caption{Modifying the Kauffman-states diagram as in Example~\ref{ex:ModifiedKSDiag}.}
\label{fig:LumpAndStabilize}
\end{figure}

\begin{example}\label{ex:ModifiedKSDiag}
Let $L$ be a marked link projection; generically we may assume that $L$ can be written as a vertical concatenation of crossings (possibly singular), maximum points, and minimum points. Such a projection $L$ is shown in the top-left corner of Figure~\ref{fig:LumpAndStabilize}.

The Kauffman-states diagram of $L$ can be viewed as being glued together from partial diagrams as shown in the top-right corner of Figure~\ref{fig:LumpAndStabilize}. We stretch the Kauffman-states diagram of $L$ as shown in the bottom-left corner of Figure~\ref{fig:LumpAndStabilize}. Then, in the spirit of Ozsv{\'a}th--Szab{\'o}'s recent work \cite{OSzNew,OSzNewer}, we perform a series of index one/two stabilizations (and handleslides of $\alpha$ circles) to obtain the diagram in the bottom-right corner of Figure~\ref{fig:LumpAndStabilize}. These stabilizations are useful when treating the Kauffman-states diagram using bordered Floer homology; they will turn out to be important in our construction as well.
\end{example}

\begin{remark}\label{rem:StabilizationPointIrrelevant}
Diagrams obtained by placing a stabilization point in a different position than the one shown in Figure~\ref{fig:LumpAndStabilize} (e.g. below a crossing as opposed to above it) are related by handleslides of $\beta$ circles.
\end{remark}

\begin{figure}
\includegraphics[scale=0.75]{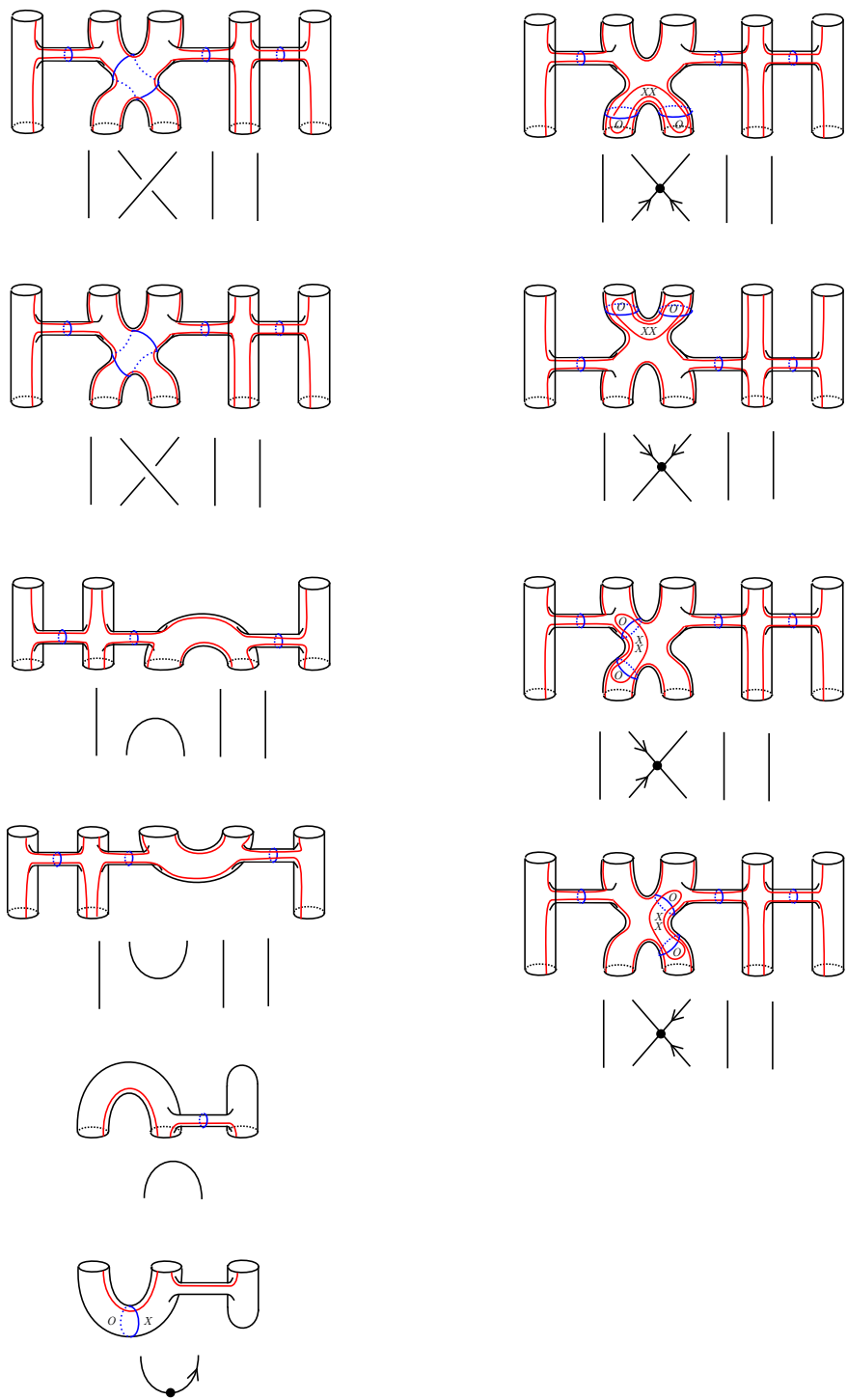}
\caption{Pieces obtained by vertically slicing the stabilized Kauffman-states diagram.}
\label{fig:KSHorizontalRows}
\end{figure}

\begin{figure}
\includegraphics[scale=0.75]{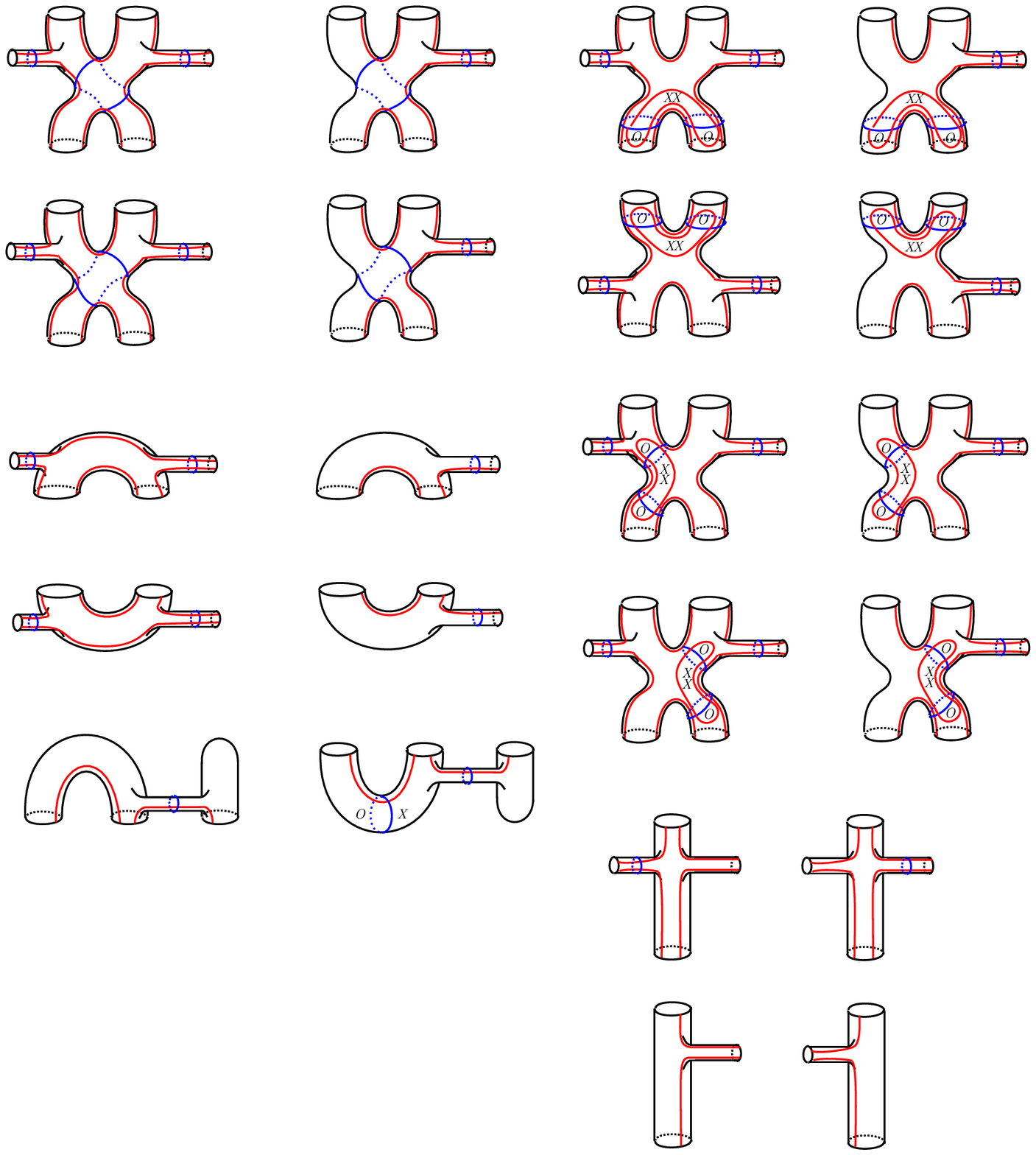}
\caption{Elementary pieces of the stabilized Kauffman-states diagram.}
\label{fig:KSDoublyLocal}
\end{figure}

\begin{example}
We can view the Kauffman-states diagram, in the version shown in the bottom-left corner of Figure~\ref{fig:LumpAndStabilize}, as glued together from the partial Heegaard diagrams shown in Figure~\ref{fig:KSHorizontalRows}. Furthermore, we can assemble each of the diagrams in Figure~\ref{fig:KSHorizontalRows} from the partial Heegaard diagrams shown in Figure~\ref{fig:KSDoublyLocal}. Ladybugs are omitted for simplicity, but in general all pieces in Figure~\ref{fig:KSDoublyLocal} are allowed to have ladybugs added.
\end{example}

\subsection{Planar pieces}

\begin{figure}
\includegraphics[scale=0.75]{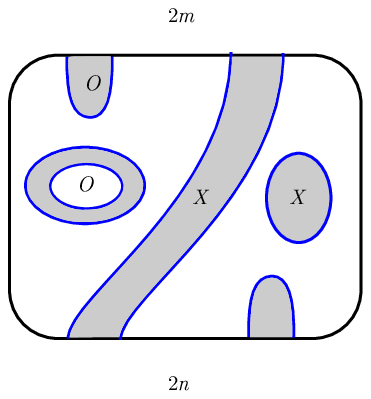}
\caption{Shading of components used in Definition~\ref{def:PlanarPiece}.}
\label{fig:PlanarPieceDef}
\end{figure}

The following definition is motivated by the structure of the planar Heegaard diagram:

\begin{definition}\label{def:PlanarPiece}
A planar piece is a partial Heegaard diagram $\Hc$ whose underlying surface is $D^2$, with $2m$ $\beta$ arcs intersecting the top part of $\partial H$ and $2n$ $\beta$ arcs intersecting the bottom part of $\partial H$, such that when the components of the complement of the $\beta$ arcs and circles in $\Hc$ are alternatingly shaded and unshaded as shown in Figure~\ref{fig:PlanarPieceDef}, there exists no path from a basepoint to $\partial \Hc$ contained entirely in an unshaded region. 
\end{definition}

\begin{example}
The partial Heegaard diagrams from Figure~\ref{fig:OrigPartialPlanarDiags}, used to construct the planar diagram in Section~\ref{sec:Planar}, are all planar pieces. 
\end{example}

\begin{figure}
\includegraphics[scale=0.75]{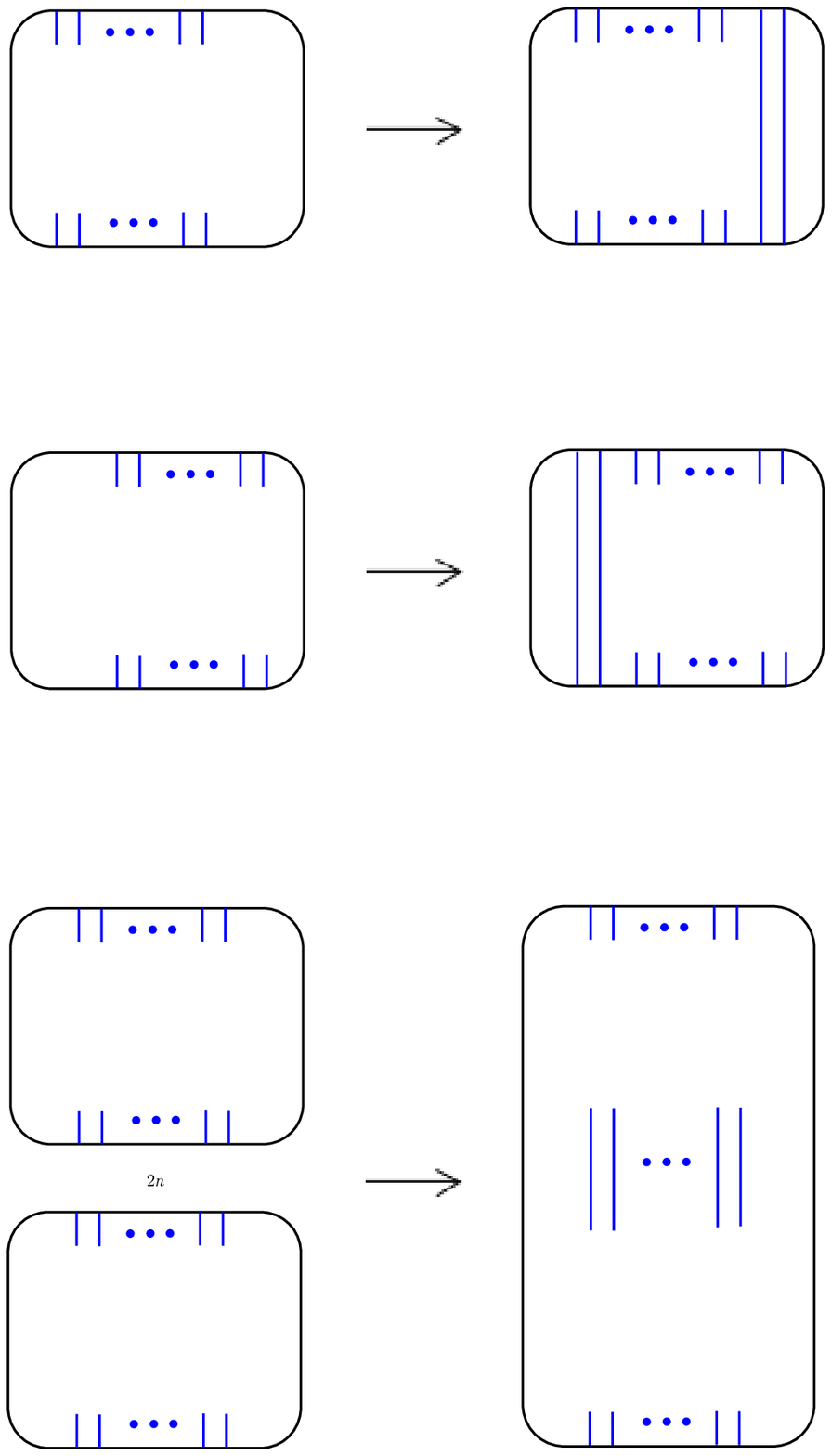}
\caption{Operations on planar pieces.}
\label{fig:PlanarPieceOps}
\end{figure}

\begin{remark}\label{rem:ModifyingPlanarPieces}
If $\Hc$ is a planar piece, then modifying $\Hc$ as shown in the top two rows Figure~\ref{fig:PlanarPieceOps} produces another planar piece. In general, one can add $2k$ parallel $\beta$ arcs to the left of $\Hc$ and $2l$ parallel $\beta$ arcs to the right of $\Hc$ and obtain another planar piece. If $\Hc_1$ and $\Hc_2$ are planar pieces with $2n$ $\beta$ arcs intersecting the top of $\Hc_1$ and the bottom of $\Hc_2$, then gluing $\Hc_1$ and $\Hc_2$ as at the bottom of Figure~\ref{fig:PlanarPieceOps} produces another planar piece. Also note that endpoints of $\beta$ arcs can be rotated two-at-a-time from the bottom to the top of a planar piece.
\end{remark}

\section{Heegaard moves local to crossings}\label{sec:Local}

The goal of this section is to transform the partial Heegaard diagrams of Figure~\ref{fig:KSDoublyLocal} into ones that can be understood using the planar pieces of Figure~\ref{fig:OrigPartialPlanarDiags} using Heegaard moves.

\begin{figure}
\includegraphics[scale=0.75]{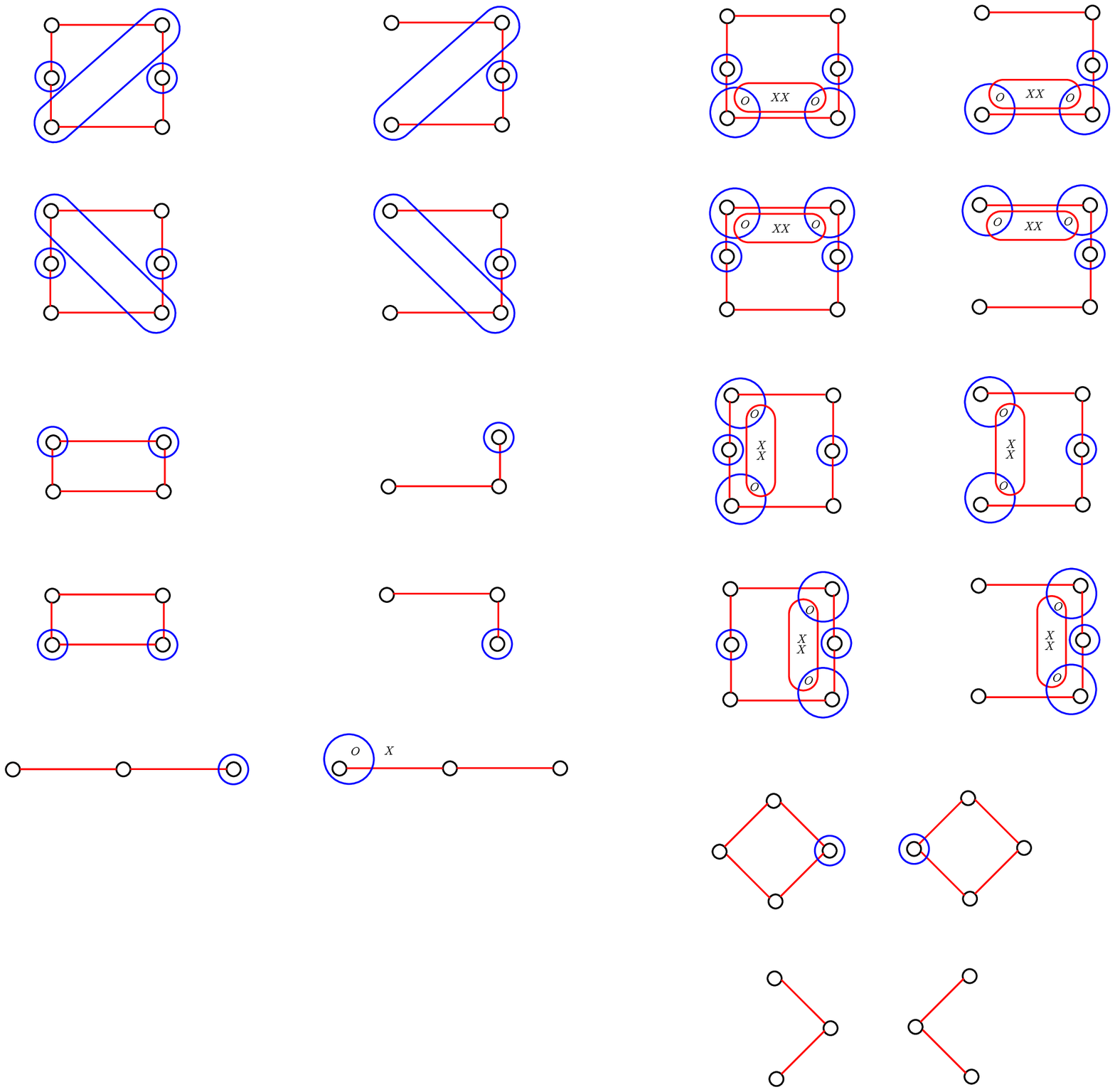}
\caption{Alternative view of the diagrams in Figure~\ref{fig:KSDoublyLocal}; each diagram is drawn on $S^2$.}
\label{fig:KSDoublyLocal2Dim}
\end{figure}

\begin{figure}
\includegraphics[scale=0.75]{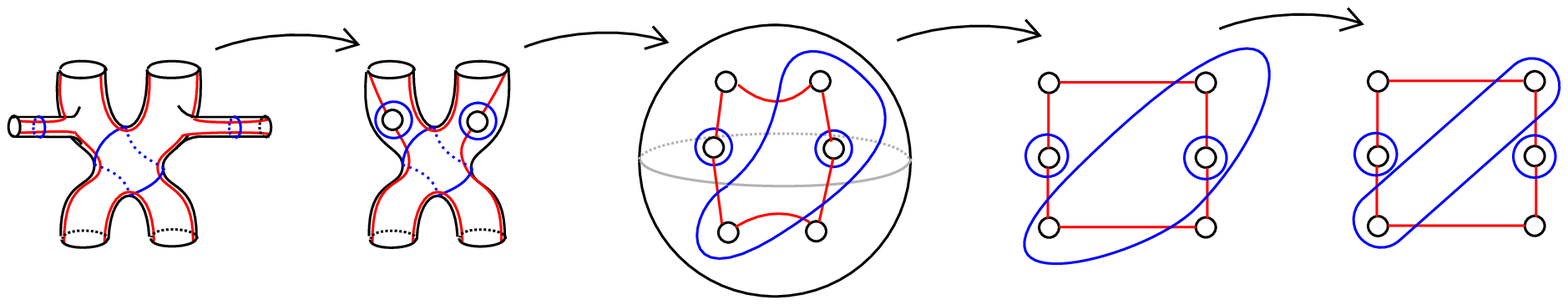}
\caption{How to go from Figure~\ref{fig:KSDoublyLocal} to Figure~\ref{fig:KSDoublyLocal2Dim}.}
\label{fig:RedrawingTheDiagrams}
\end{figure}

To describe these moves, it is helpful to re-draw the diagrams of Figure~\ref{fig:KSDoublyLocal} in a more planar way. Figure~\ref{fig:KSDoublyLocal2Dim} shows an alternative view of the diagrams of Figure~\ref{fig:KSDoublyLocal} (with a few $\beta$ handleslides done to make the diagrams look nicer; these come from moving the stabilization points as in Remark~\ref{rem:StabilizationPointIrrelevant}). All diagrams are assumed to be drawn on $S^2$; we will not draw these spheres to save space. Figure~\ref{fig:RedrawingTheDiagrams} shows how to go from one view to the other in one example.

\begin{figure}
\includegraphics[scale=0.75]{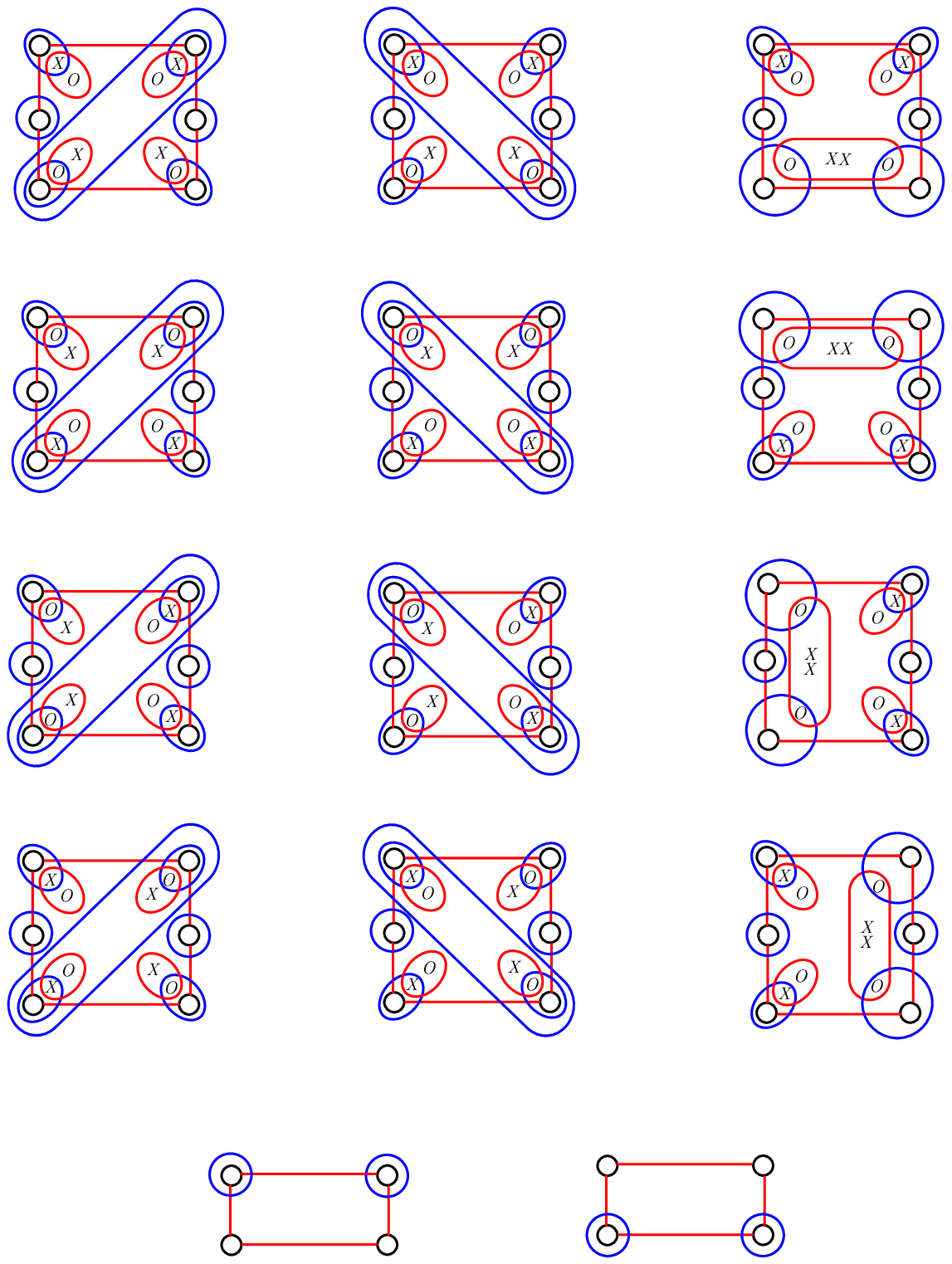}
\caption{Local diagrams that are the starting point of Construction~\ref{constr:Local} (generic case). Several have ladybugs added when compared with Figure~\ref{fig:KSDoublyLocal2Dim}.}
\label{fig:KSDLWithLadybugs}
\end{figure}

\begin{figure}
\includegraphics[scale=0.75]{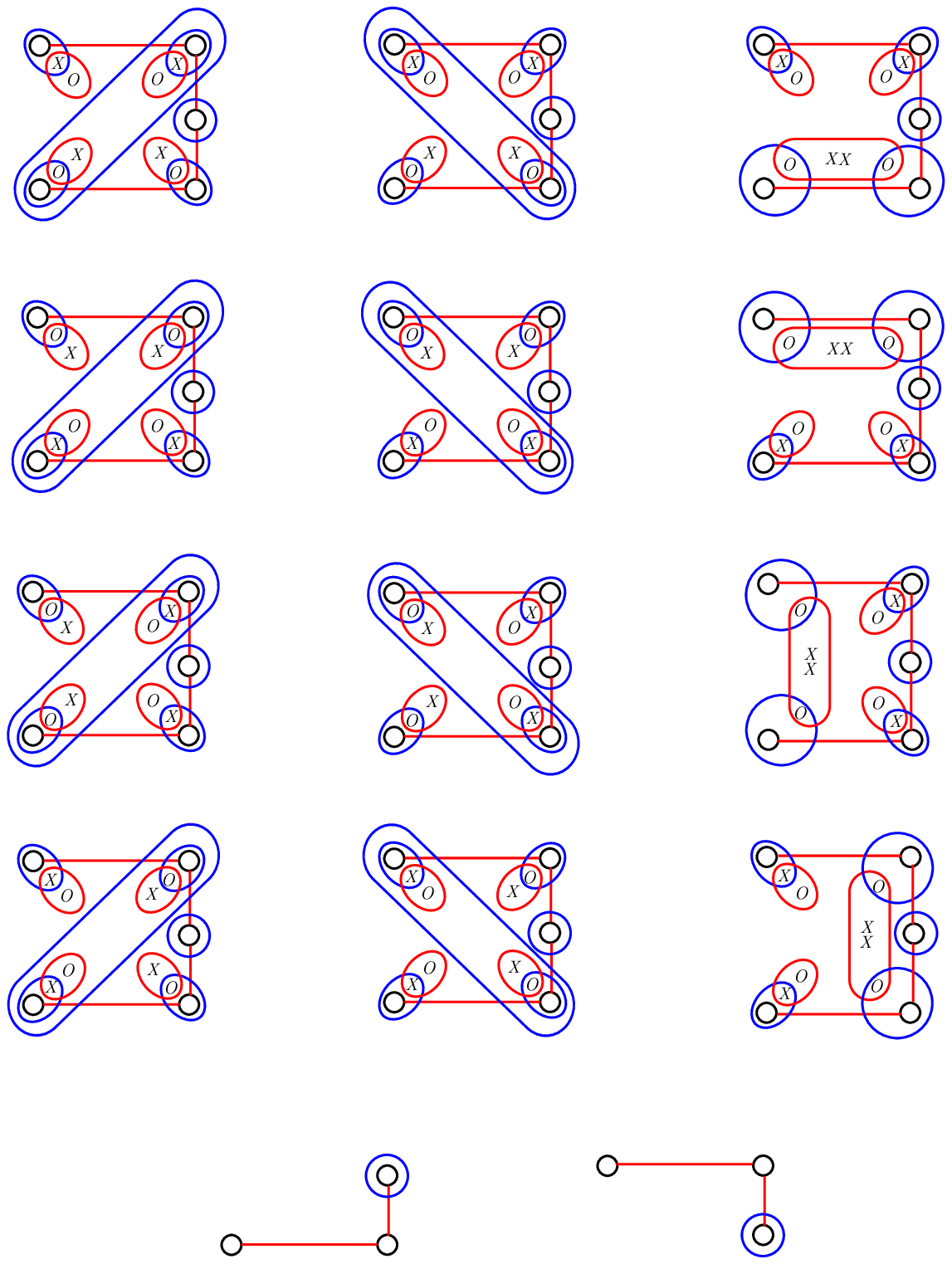}
\caption{Local diagrams that are the starting point of Construction~\ref{constr:Local} (non-generic case).}
\label{fig:NonGenericParts}
\end{figure}

We now add the ladybugs that we will need (these were assumed to be potentially present the whole time). Some of the diagrams from Figure~\ref{fig:KSDoublyLocal2Dim} that get ladybugs are shown (along with two diagrams that get no ladybugs) in Figure~\ref{fig:KSDLWithLadybugs}. Nonsingular crossings can have four possible orientation patterns, and the $X$ and $O$ basepoints in the ladybugs depend on the orientation.

The rest of the diagrams that get ladybugs are shown in Figure~\ref{fig:NonGenericParts}, again with two diagrams that get no ladybugs. The diagrams in Figure~\ref{fig:KSDLWithLadybugs} and Figure~\ref{fig:NonGenericParts} are the ones we will transform using Construction~\ref{constr:Local}. 

\begin{figure}
\includegraphics[scale=0.75]{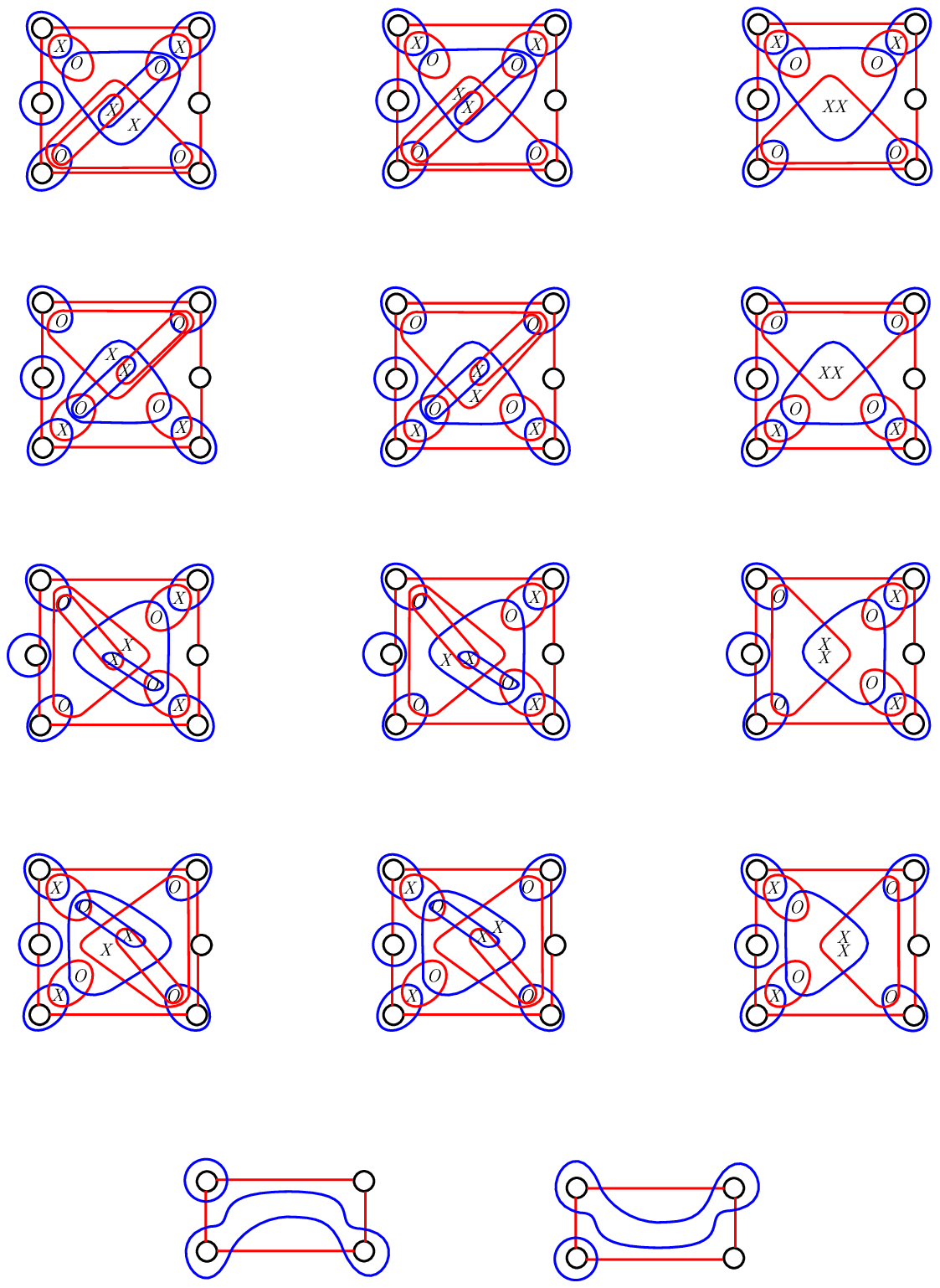}
\caption{Local diagrams that are the end result of Construction~\ref{constr:Local} (generic case).}
\label{fig:TransformedLocalParts}
\end{figure}

\begin{figure}
\includegraphics[scale=0.75]{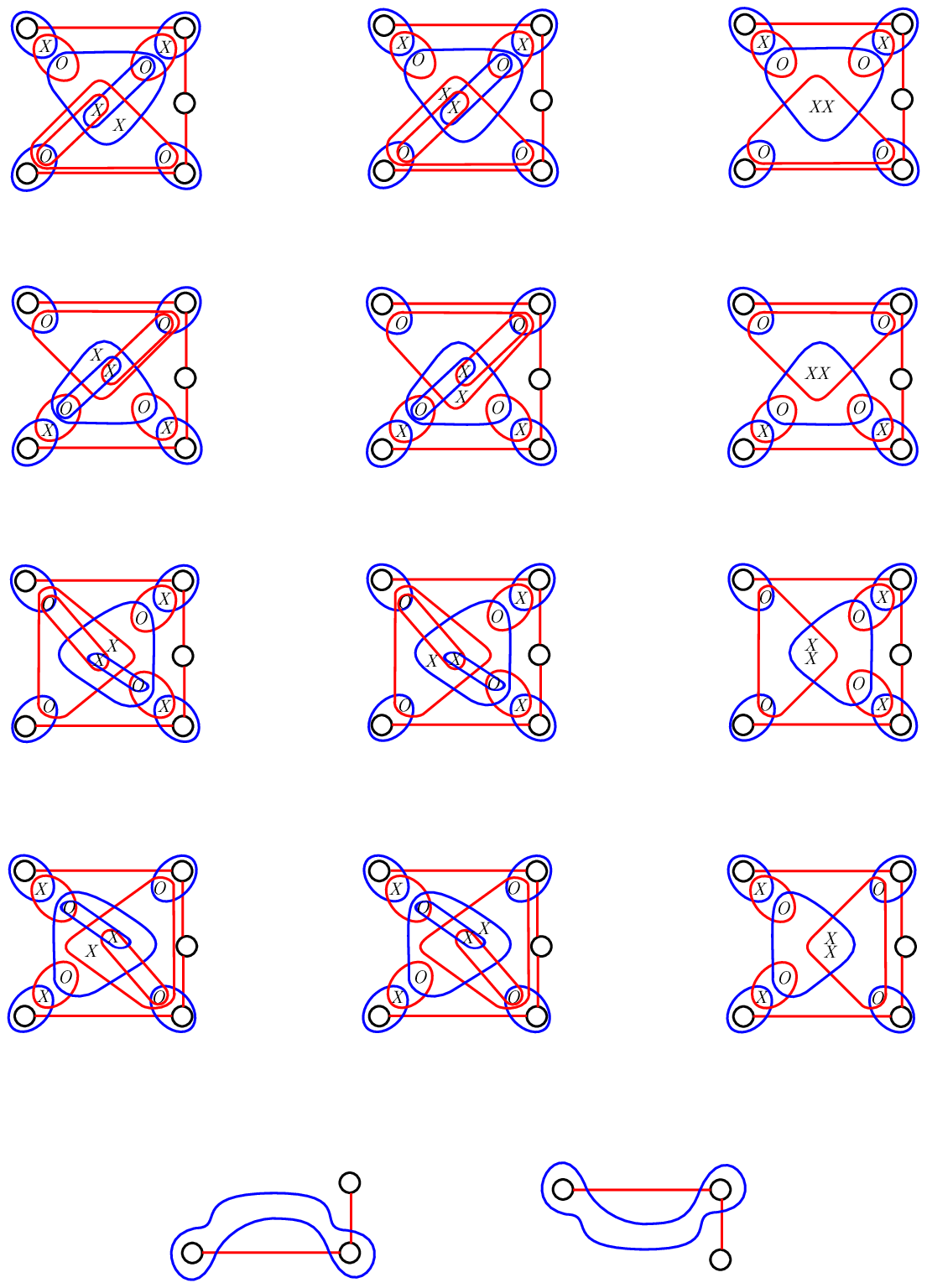}
\caption{Local diagrams that are the end result of Construction~\ref{constr:Local}: (non-generic case).}
\label{fig:NonGenericTransformed}
\end{figure}


\begin{figure}
\includegraphics[scale=0.75]{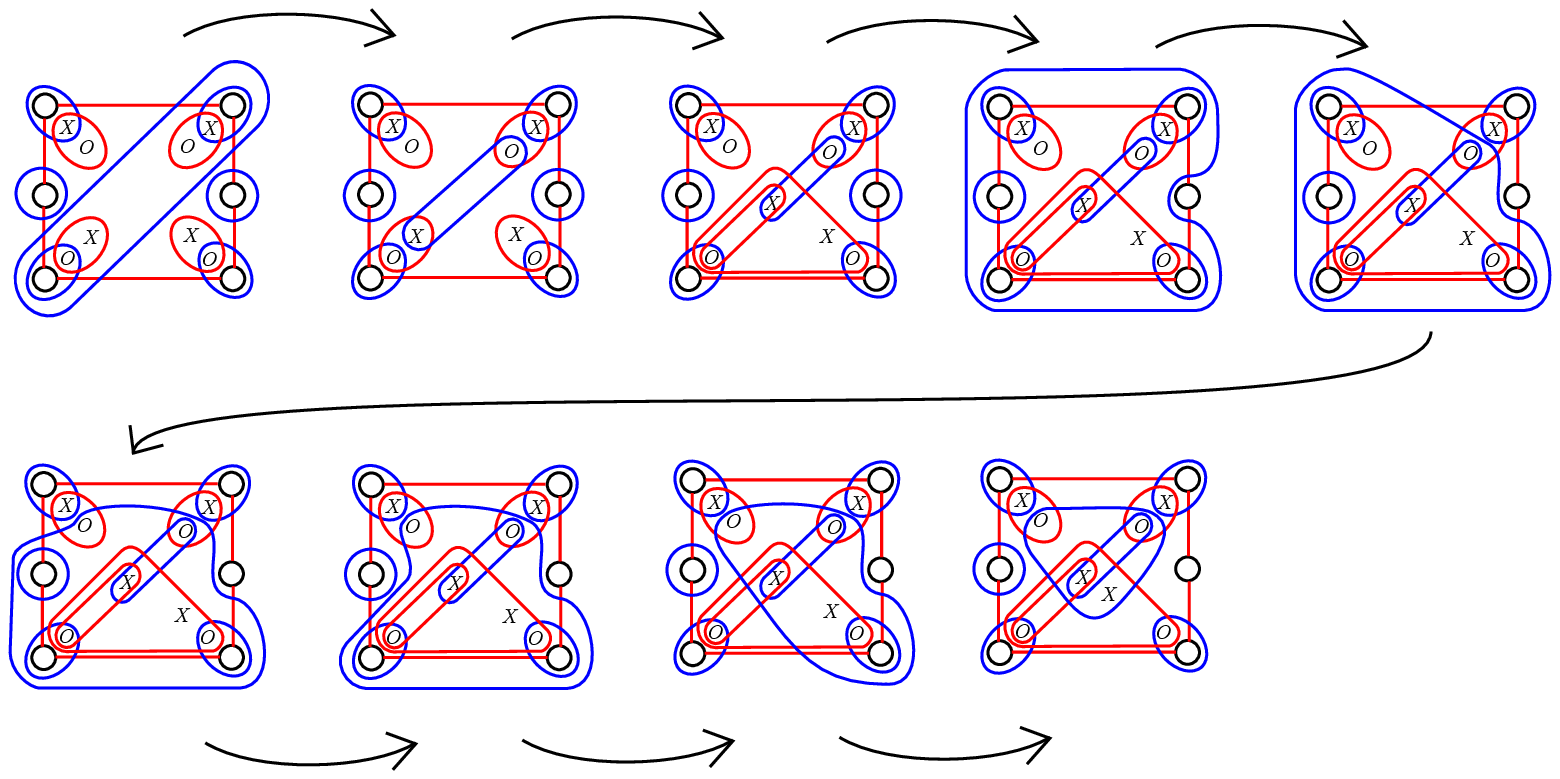}
\caption{Heegaard moves used in Construction~\ref{constr:Local}.}
\label{fig:LocalConstruction1}
\end{figure}

\begin{figure}
\includegraphics[scale=0.75]{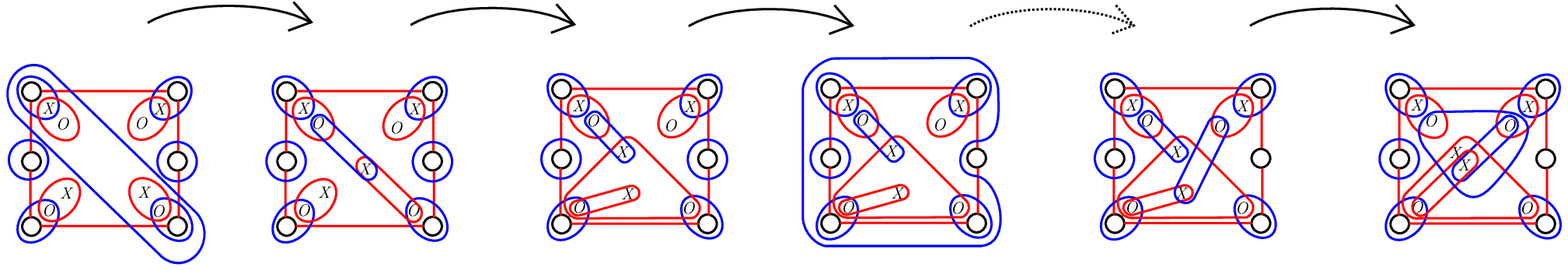}
\caption{Heegaard moves used in Construction~\ref{constr:Local}.}
\label{fig:LocalConstruction3}
\end{figure}

\begin{figure}
\includegraphics[scale=0.75]{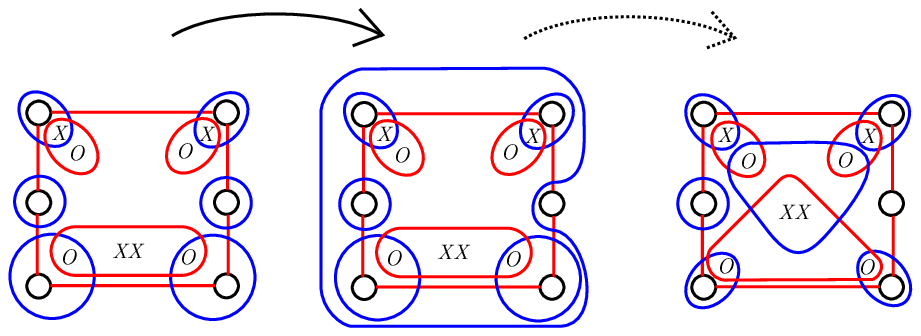}
\caption{Heegaard moves used in Construction~\ref{constr:Local}.}
\label{fig:LocalConstruction5}
\end{figure}


\begin{figure}
\includegraphics[scale=0.75]{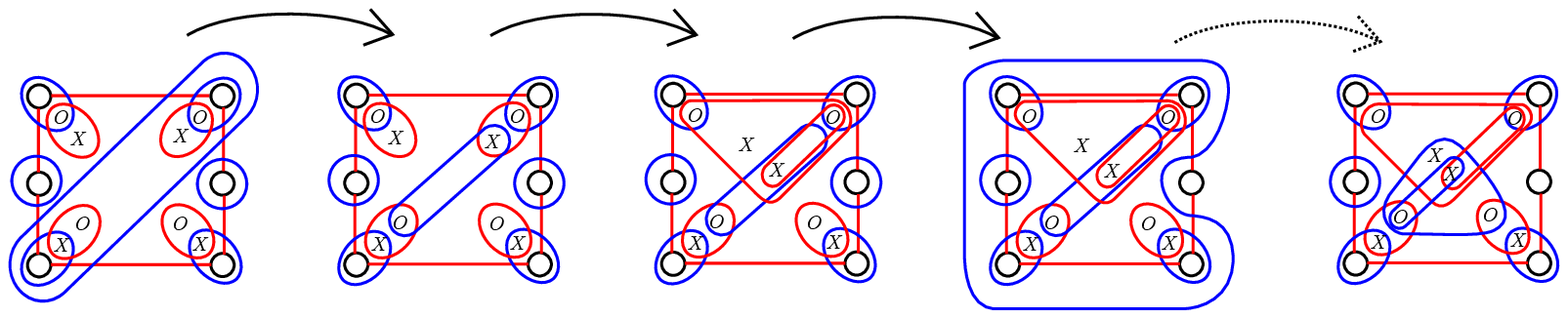}
\caption{Heegaard moves used in Construction~\ref{constr:Local}.}
\label{fig:LocalConstruction5A}
\end{figure}

\begin{figure}
\includegraphics[scale=0.75]{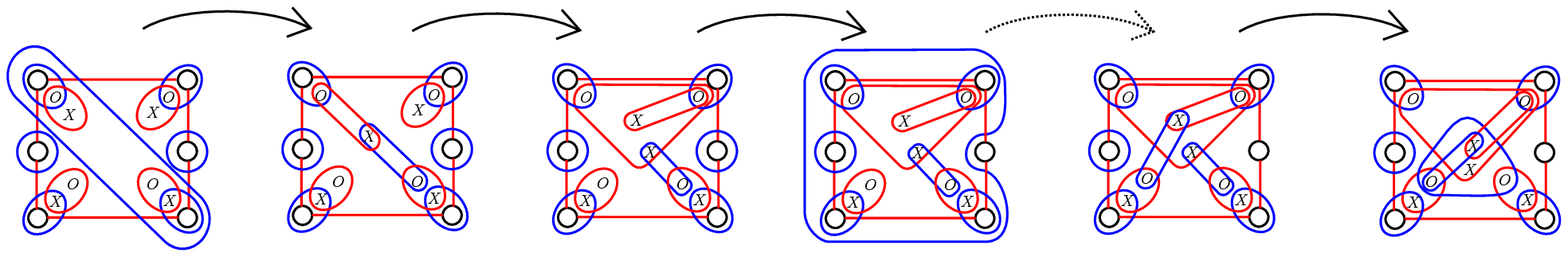}
\caption{Heegaard moves used in Construction~\ref{constr:Local}.}
\label{fig:LocalConstruction5B}
\end{figure}

\begin{figure}
\includegraphics[scale=0.75]{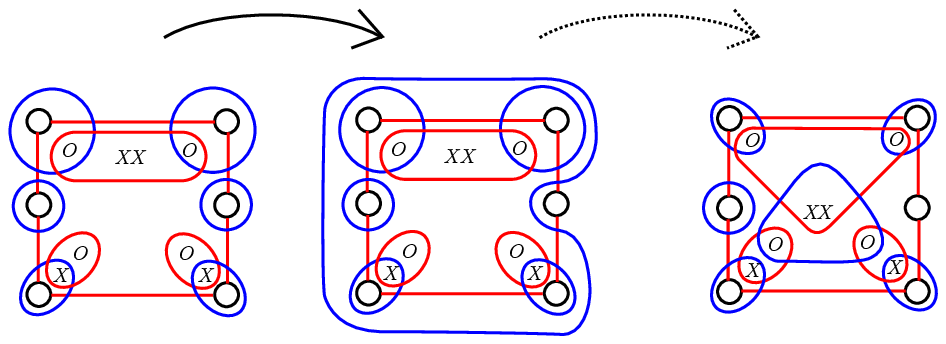}
\caption{Heegaard moves used in Construction~\ref{constr:Local}.}
\label{fig:LocalConstruction6}
\end{figure}


\begin{figure}
\includegraphics[scale=0.75]{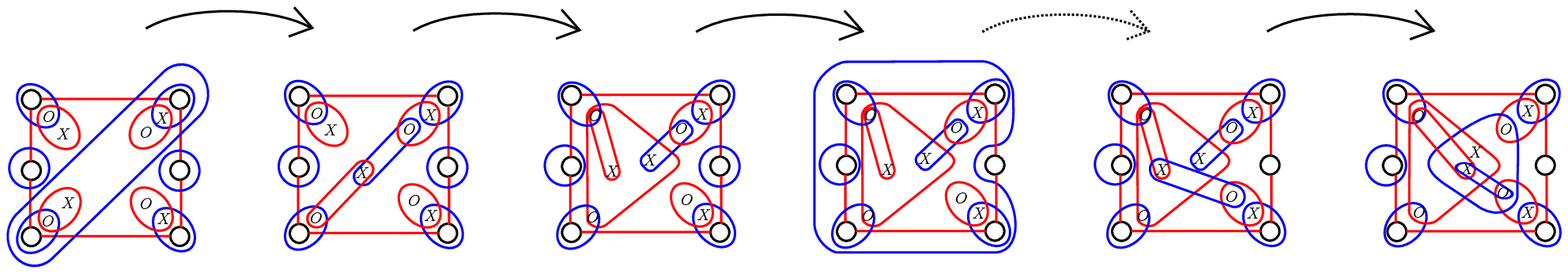}
\caption{Heegaard moves used in Construction~\ref{constr:Local}.}
\label{fig:LocalConstruction6A}
\end{figure}

\begin{figure}
\includegraphics[scale=0.75]{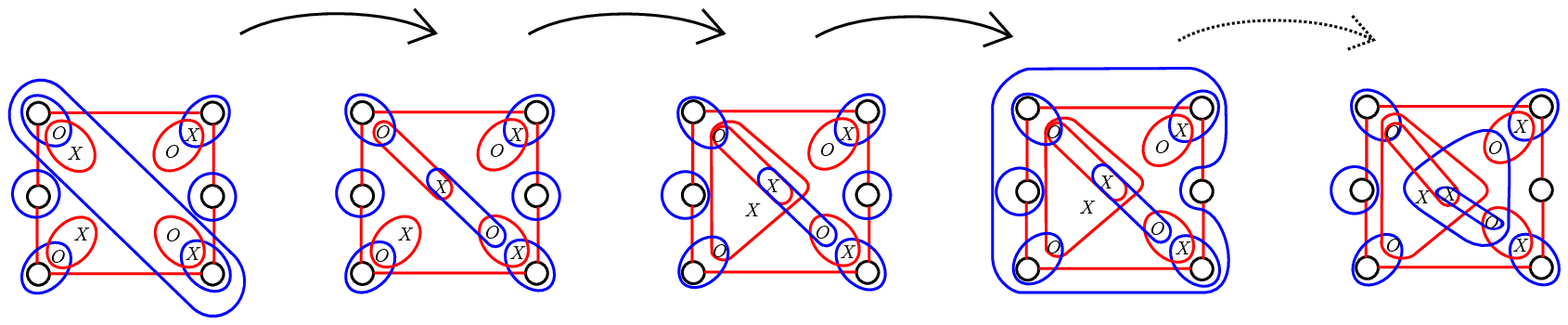}
\caption{Heegaard moves used in Construction~\ref{constr:Local}.}
\label{fig:LocalConstruction6B}
\end{figure}

\begin{figure}
\includegraphics[scale=0.75]{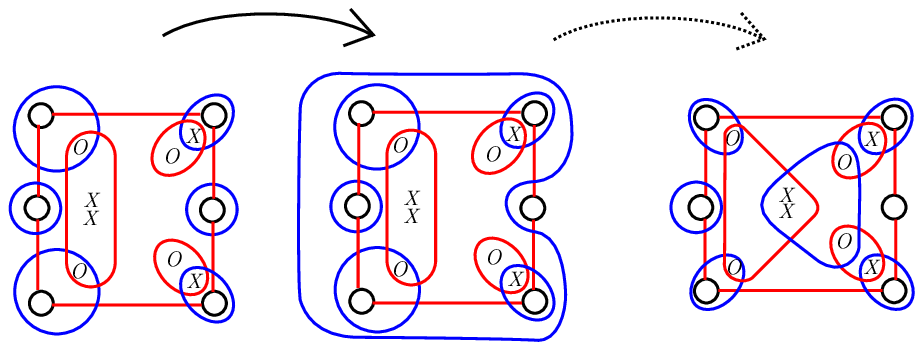}
\caption{Heegaard moves used in Construction~\ref{constr:Local}.}
\label{fig:LocalConstruction7}
\end{figure}


\begin{figure}
\includegraphics[scale=0.75]{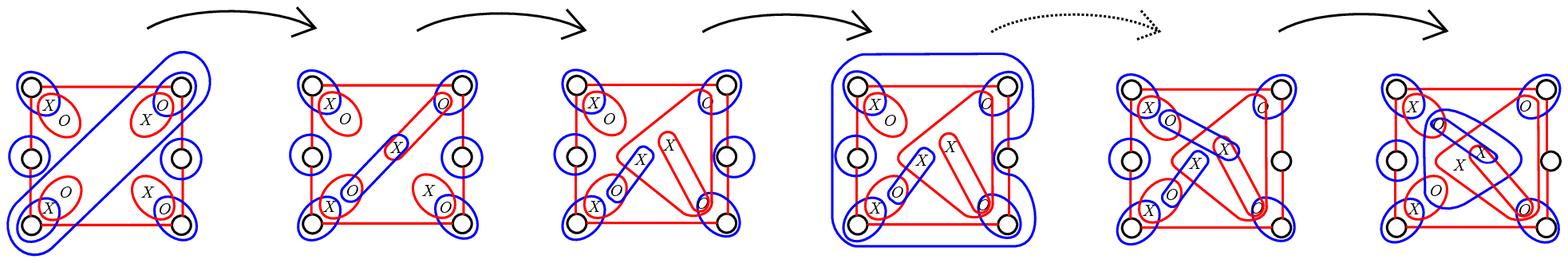}
\caption{Heegaard moves used in Construction~\ref{constr:Local}.}
\label{fig:LocalConstruction7A}
\end{figure}

\begin{figure}
\includegraphics[scale=0.75]{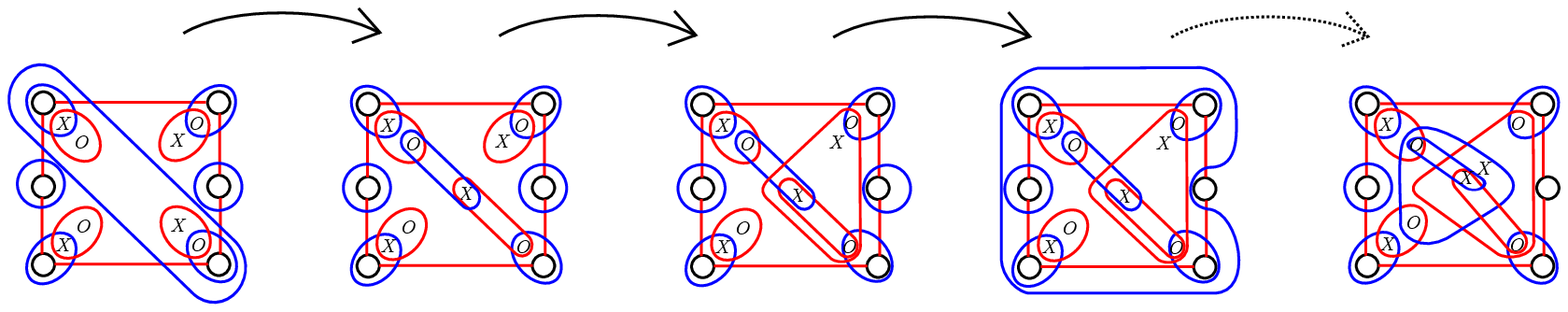}
\caption{Heegaard moves used in Construction~\ref{constr:Local}.}
\label{fig:LocalConstruction7B}
\end{figure}

\begin{figure}
\includegraphics[scale=0.75]{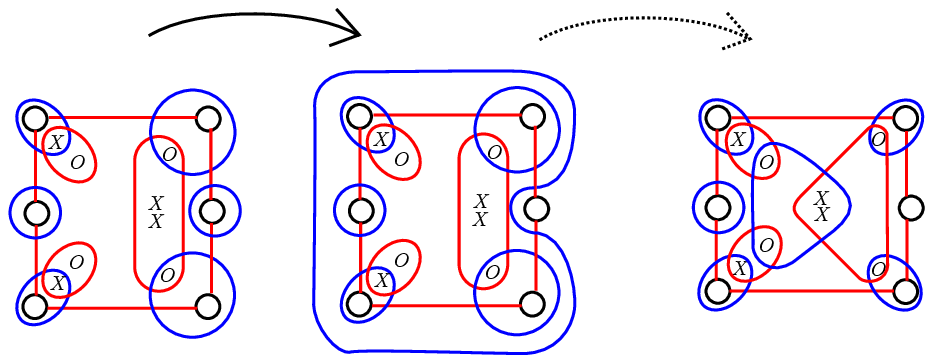}
\caption{Heegaard moves used in Construction~\ref{constr:Local}.}
\label{fig:LocalConstruction8}
\end{figure}


\begin{figure}
\includegraphics[scale=0.75]{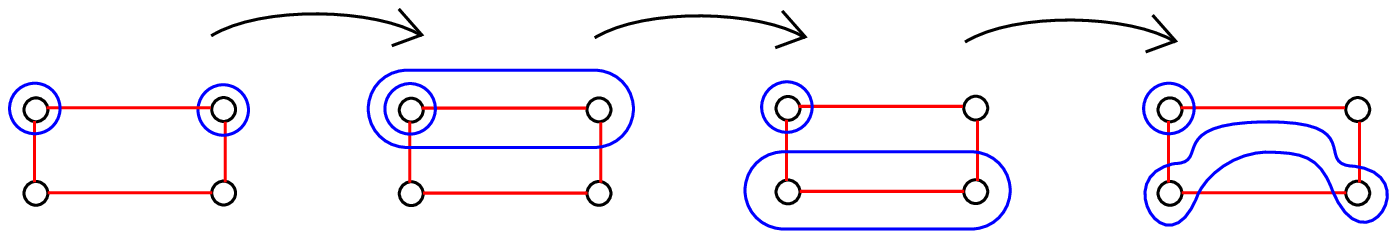}
\caption{Heegaard moves used in Construction~\ref{constr:Local}.}
\label{fig:LocalConstruction9}
\end{figure}

\begin{figure}
\includegraphics[scale=0.75]{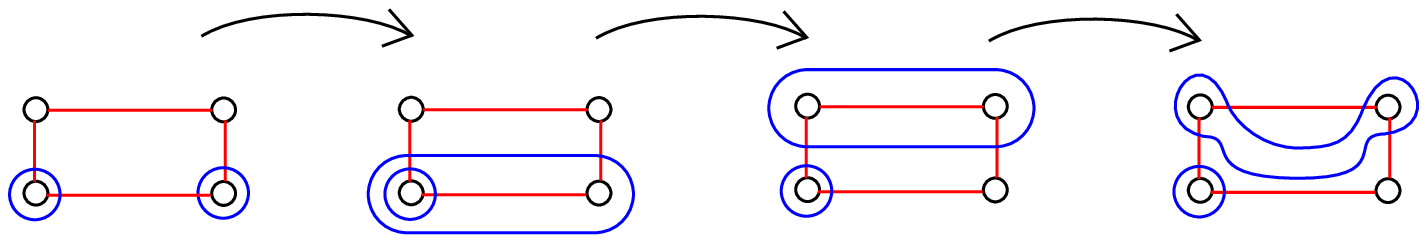}
\caption{Heegaard moves used in Construction~\ref{constr:Local}.}
\label{fig:LocalConstruction10}
\end{figure}

The following construction transforms the diagrams of Figure~\ref{fig:KSDLWithLadybugs} into the diagrams of Figure~\ref{fig:TransformedLocalParts}. The same steps apply to transform the diagrams of Figure~\ref{fig:NonGenericParts} into the diagrams of Figure~\ref{fig:NonGenericTransformed}.

\begin{construction}\label{constr:Local}
Starting from the diagram at the left of Figure~\ref{fig:LocalConstruction1} (row 1, column 1 of Figure~\ref{fig:KSDLWithLadybugs}), handleslide the ``diagonal'' $\beta$ circle to get the second diagram of Figure~\ref{fig:LocalConstruction1}. Handleslide an $\alpha$ circle to get the third diagram of Figure~\ref{fig:LocalConstruction1}. Wrap the rightmost $\beta$ circle around the back of $S^2$ to get the fourth diagram of Figure~\ref{fig:LocalConstruction1}. Finally, perform several handleslides of this $\beta$ circle, as shown in the subsequent diagrams of Figure~\ref{fig:LocalConstruction1}. 

The rest of the diagrams of Figure~\ref{fig:KSDLWithLadybugs} are analyzed similarly; they are shown in Figures~\ref{fig:LocalConstruction3}--\ref{fig:LocalConstruction10}. To save space, long sequences of handleslides of the same $\beta$ circle like the one occurring in Figure~\ref{fig:LocalConstruction1} are omitted; we indicate them with dotted arrows.
\end{construction}

\begin{figure}
\includegraphics[scale=0.75]{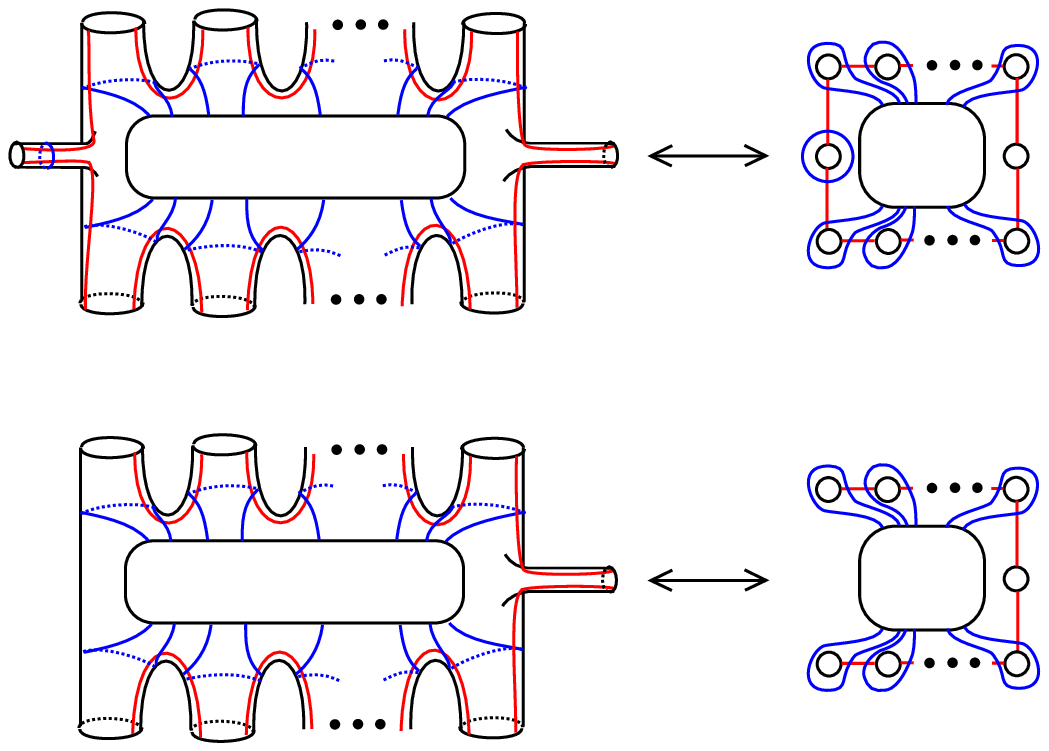}
\caption{Standard frames for inserting planar pieces.}
\label{fig:StandardFrame}
\end{figure}

\begin{figure}
\includegraphics[scale=0.75]{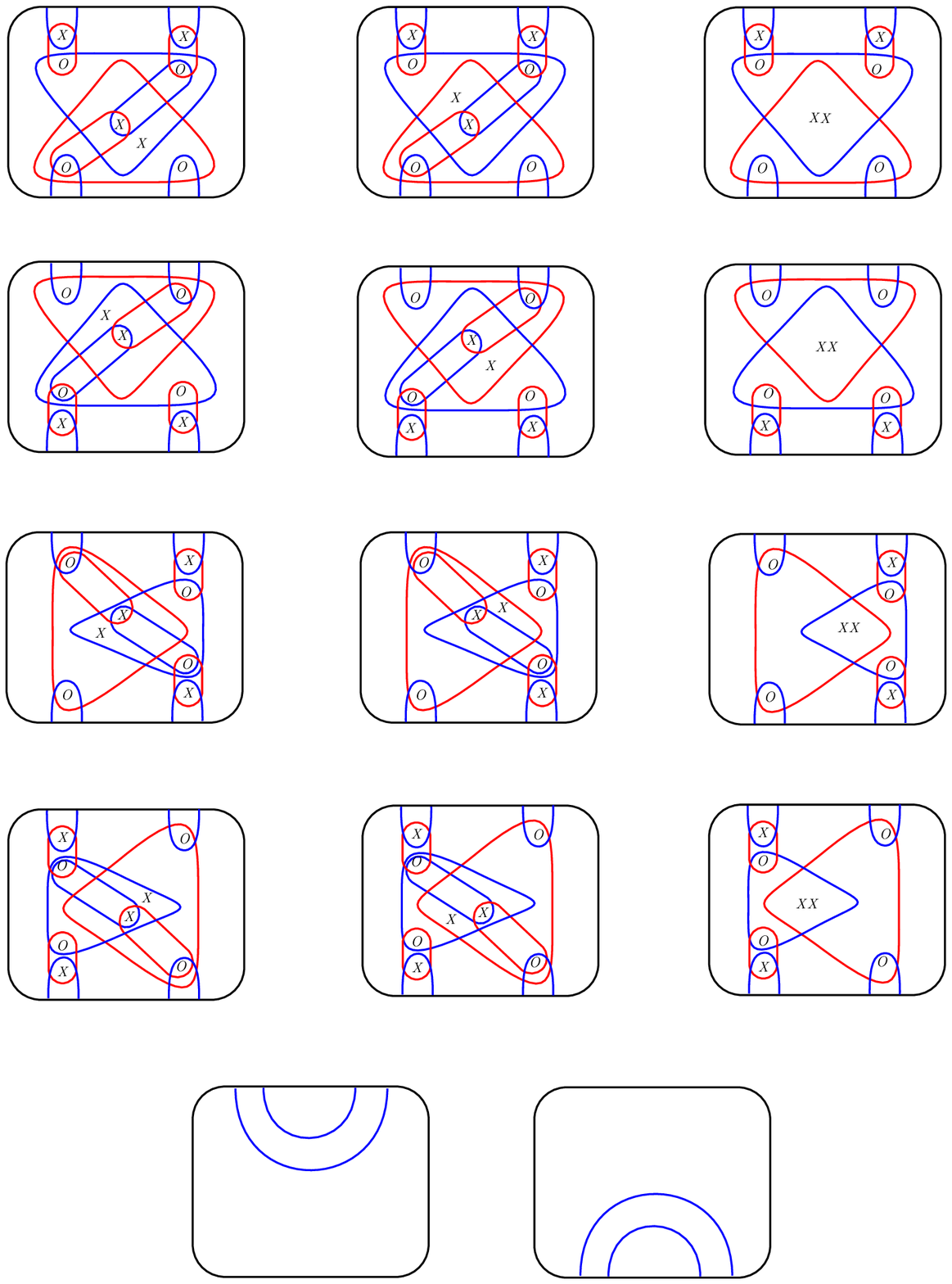}
\caption{Pieces of the planar diagram to glue to the diagrams in Figure~\ref{fig:StandardFrame}.}
\label{fig:LotsOfPlanarPieces}
\end{figure}

The diagrams in Figure~\ref{fig:TransformedLocalParts} all have a special form. If $\Hc_P$ is a planar piece with $2m$ endpoints of $\beta$ arcs along the top of $\Hc_{\planar}$ and $2n$ endpoints of $\beta$ arcs along the bottom, we can form a new partial Heegaard diagram by gluing $\Hc_P$ into the diagram shown in the top row of Figure~\ref{fig:StandardFrame}. We get the diagrams in Figure~\ref{fig:NonGenericTransformed} by gluing $\Hc_P$ into the diagram shown on in the bottom row of Figure~\ref{fig:StandardFrame} instead. All the diagrams in Figure~\ref{fig:TransformedLocalParts} can be obtained by applying this procedure to a planar piece from Figure~\ref{fig:LotsOfPlanarPieces}; these pieces are the possible local pieces of the planar diagram when one does not mod out by rotation.

\begin{figure}
\includegraphics[scale=0.75]{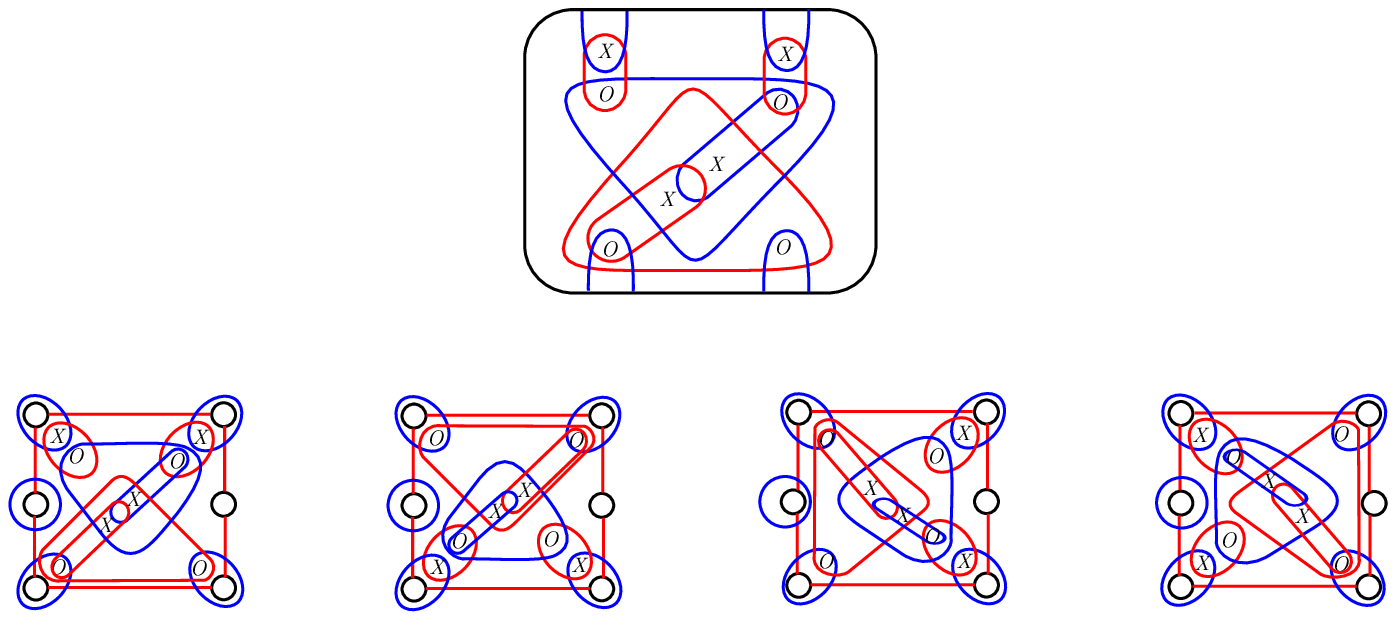}
\caption{The local piece of the planar diagram for a resolved crossing, and the different ways of gluing it to the top diagram of Figure~\ref{fig:StandardFrame}.}
\label{fig:ResolvedCrossingPlanarIntro}
\end{figure}

\begin{figure}
\includegraphics[scale=0.75]{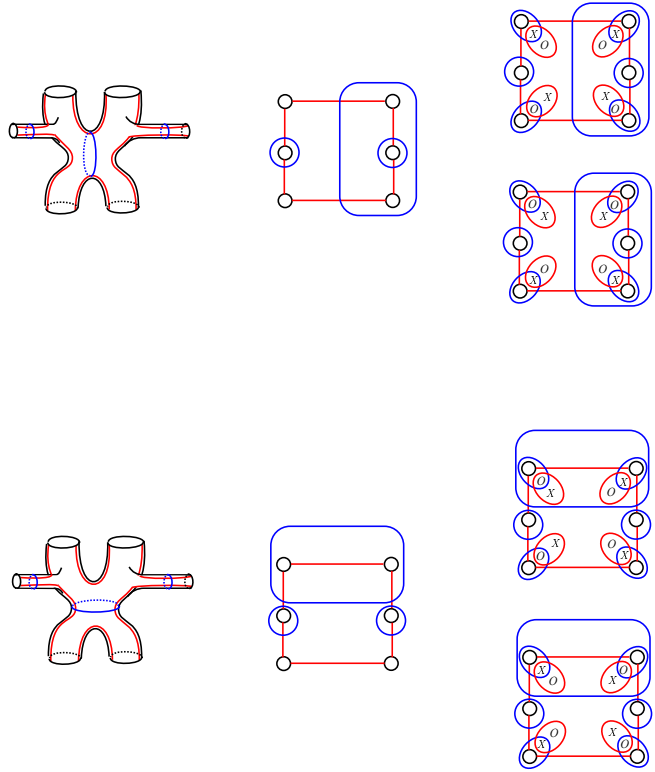}
\caption{Local Kauffman-states diagrams for a resolved crossing, with and without ladybugs.}
\label{fig:ResolvedCrossingKSIntro}
\end{figure}

\begin{figure}
\includegraphics[scale=0.75]{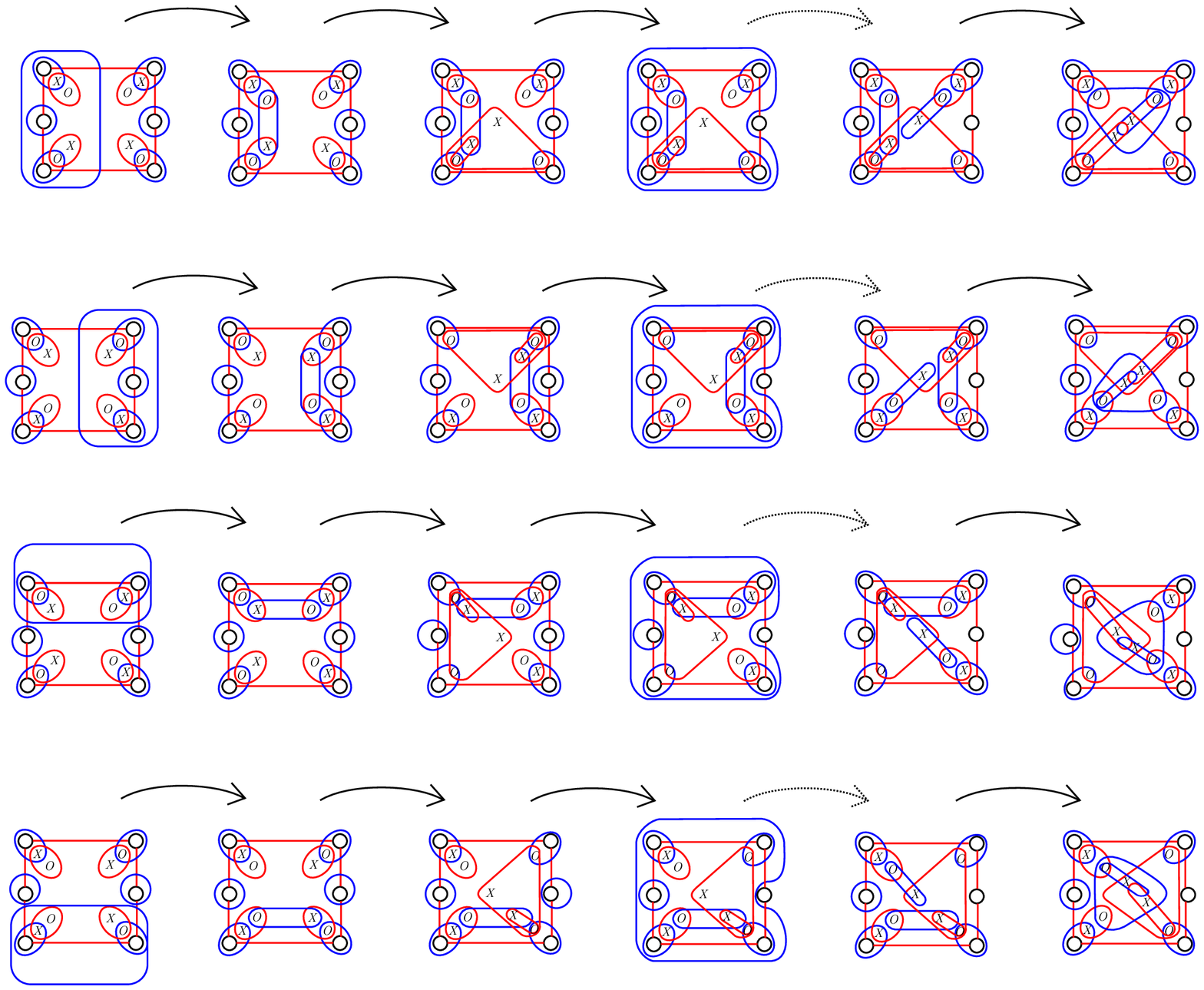}
\caption{Analogue of Construction~\ref{constr:Local} for resolved crossings.}
\label{fig:ResolvedCrossing}
\end{figure}

\begin{remark}
When working with the cube of resolutions, there is a third variant of the planar diagram that appears for resolutions of crossings. This diagram is shown in Figure~\ref{fig:ResolvedCrossingPlanarIntro}, along with the four ways to insert it into the top diagram of Figure~\ref{fig:StandardFrame} depending on orientations. Construction~\ref{constr:Local} can be modified to obtain these diagram from the Kauffman-states diagrams shown on the right of Figure~\ref{fig:ResolvedCrossingKSIntro}. We show the required Heegaard moves in Figure~\ref{fig:ResolvedCrossing}; they are very similar to the ones appearing in Construction~\ref{constr:Local}.
\end{remark}

\section{Global diagrams}\label{sec:Global}

\begin{figure}
\includegraphics[scale=0.75]{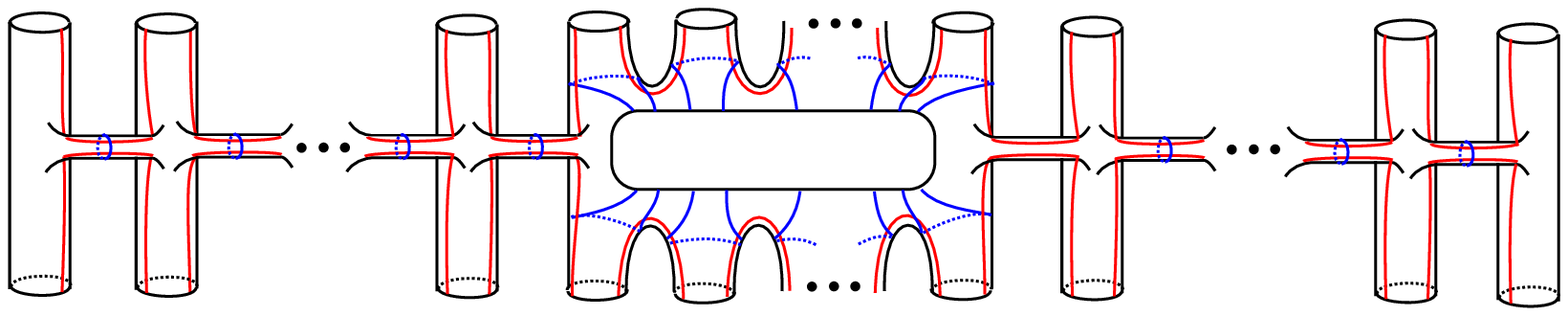}
\caption{The Heegaard diagram $\Hc_{\fr}$ used in Construction~\ref{constr:HorizExt}.}
\label{fig:ExtendedFrame}
\end{figure}

\begin{figure}
\includegraphics[scale=0.75]{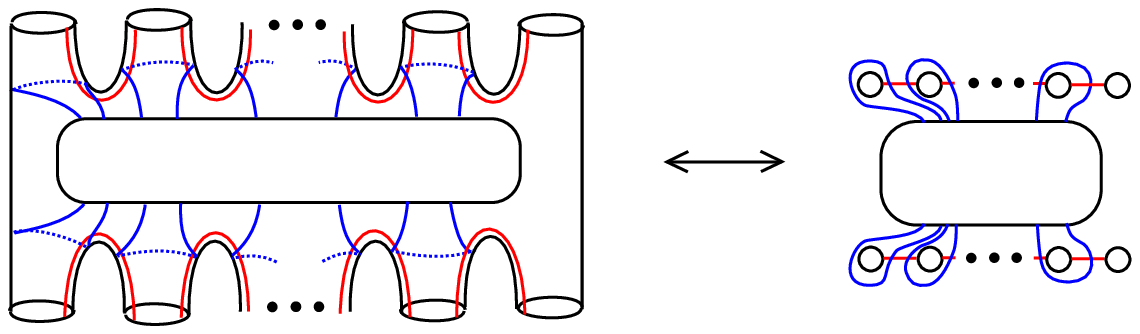}
\caption{The Heegaard diagram $\Hc'_{\fr}$ used in Construction~\ref{constr:HorizExt}.}
\label{fig:ClosedFrame}
\end{figure}

Let $\Hc_P$ be a planar piece with $2m$ top endpoints and $2n$ bottom endpoints. Let $\Hc_F$ be the partial Heegaard diagram of Figure~\ref{fig:ExtendedFrame}, with endpoints chosen so that it makes sense to glue $\Hc_F$ into $\Hc_F$. Let $\Hc'_P$ denote the planar piece obtained by inserting $2k \geq 0$ parallel $\beta$ arcs at the left of $\Hc_{\planar}$ and $2l \geq 2$ parallel $\beta$ arcs at the right of $\Hc_P$ as in Remark~\ref{rem:ModifyingPlanarPieces}. Finally, let $\Hc'_P$ denote the Heegaard diagram shown in Figure~\ref{fig:ClosedFrame} with $2(k+m+l)$ endpoints of $\beta$ arcs at the top and $2(k+n+l)$ endpoints of $\beta$ arcs at the bottom. We will give a sequence of Heegaard moves transforming $\Hc_F \cup \Hc_P$ into $\Hc'_F \cup \Hc'_P$. 

\begin{figure}
\includegraphics[scale=0.75]{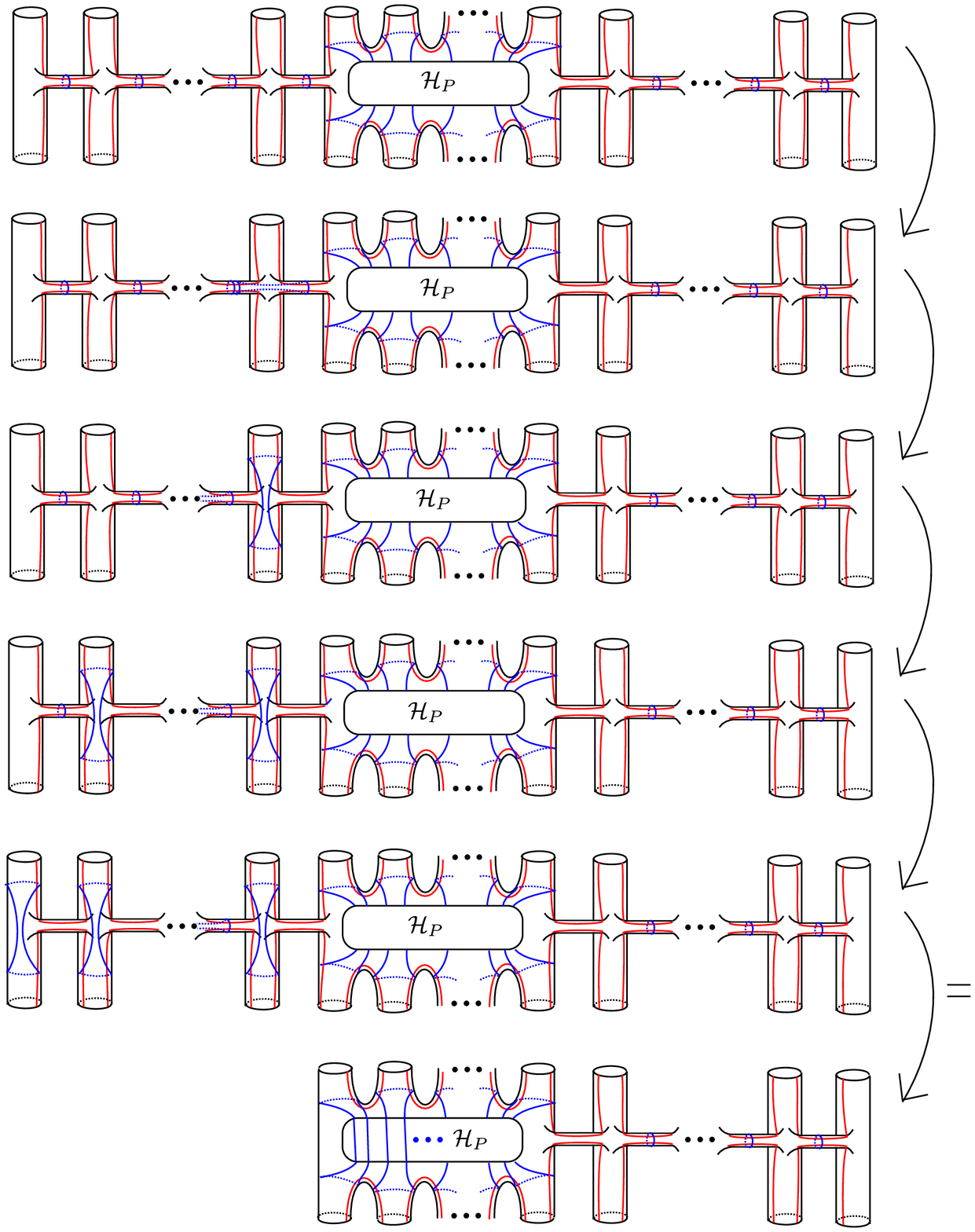}
\caption{The first part of Construction~\ref{constr:HorizExt}.}
\label{fig:HorizConstGen}
\end{figure}

\begin{figure}
\includegraphics[scale=0.75]{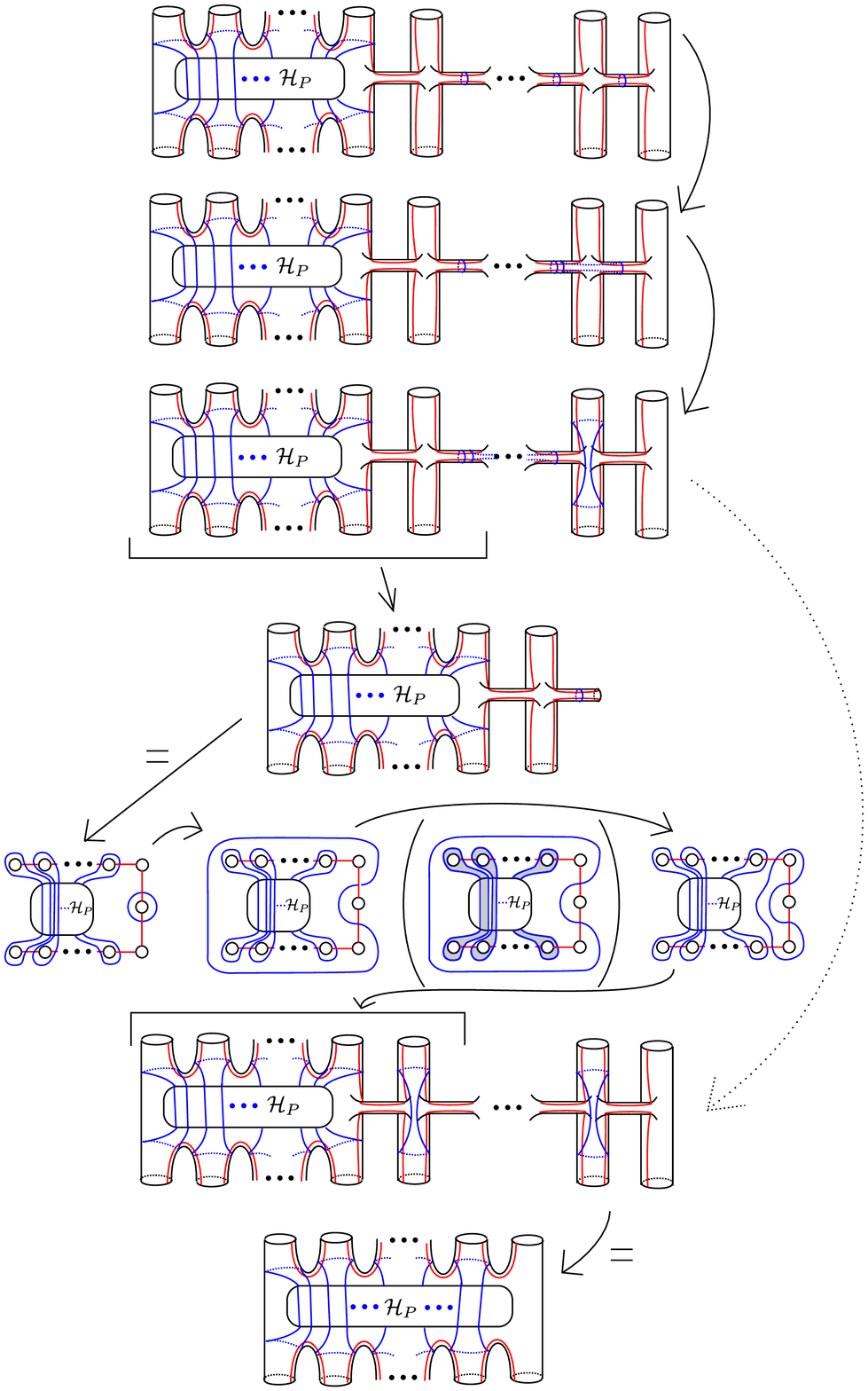}
\caption{The second part of Construction~\ref{constr:HorizExt}.}
\label{fig:HorizConstGenPart2}
\end{figure}

\begin{construction}\label{constr:HorizExt}
Starting from $\Hc_F \cup \Hc_P$ as shown in the top row of Figure~\ref{fig:HorizConstGen}, perform a sequence of handleslides on the $\beta$ curves on the left side of the diagram as shown in the subsequent rows of Figure~\ref{fig:HorizConstGen}. Next, do the same with the $\beta$ curves on the right side of the diagram, as shown in Figure~\ref{fig:HorizConstGenPart2}.

From the diagram in the fourth row of Figure~\ref{fig:HorizConstGenPart2}, wrap a $\beta$ circle around the back of $S^2$, then handleslide it through $\Hc_P$ as in the fifth row of Figure~\ref{fig:HorizConstGenPart2}. The handleslide through $\Hc_P$ is the crucial step; it is possible because $\Hc_P$ was assumed to be a planar piece. The resulting Heegaard diagram, shown at the bottom of Figure~\ref{fig:HorizConstGenPart2}, is $\Hc'_F \cup \Hc'_P$.
\end{construction}

\begin{figure}
\includegraphics[scale=0.75]{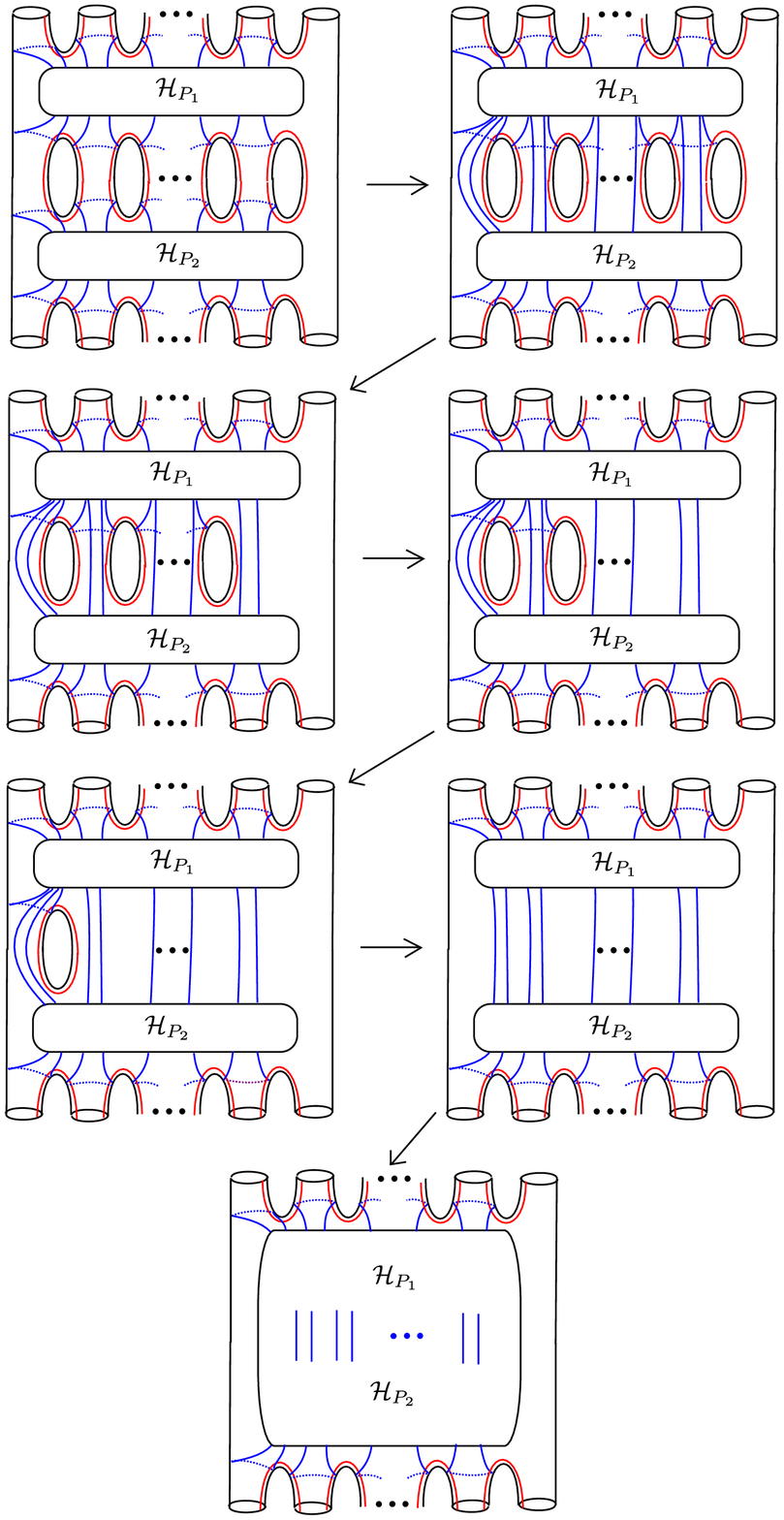}
\caption{Construction~\ref{constr:VertExt}.}
\label{fig:VertConst}
\end{figure}

Next, let $\Hc_{P_1}$ and $\Hc_{P_2}$ be planar pieces such that it makes sense to glue $\Hc_{P_1}$ on top of $\Hc_{P_2}$ as in Remark~\ref{rem:ModifyingPlanarPieces}. Let $\Hc_P$ denote the planar piece obtained by gluing $\Hc_{P_1}$ and $\Hc_{P_2}$ like this. We can glue $\Hc_P$ into the diagram $\Hc'_{F}$ of Figure~\ref{fig:ClosedFrame} above; we can also glue $\Hc_{P_i}$ into their own instances of $\Hc'_F$ and then glue the results as in the top row of Figure~\ref{fig:VertConst}. The next construction gives a sequence of Heegaard moves relating these two diagrams.

\begin{construction}\label{constr:VertExt}
Starting from the top row of Figure~\ref{fig:VertConst}, start by handlesliding $\beta$ circles to get the second row of Figure~\ref{fig:VertConst}. Then, starting at the right side of the diagram, perform a sequence of index one/two destabilizations. At the end of this procedure, one reaches $\Hc_P \cup \Hc'_F$, shown at the bottom of Figure~\ref{fig:VertConst}.
\end{construction}

We now give our main construction, illustrated by the trefoil example we have been considering. Let $L$ be a link projection, as in Example~\ref{ex:ModifiedKSDiag}, assume that $L$ can be written as a vertical concatenation of crossings (possibly singular), maximum points, and minimum points (after a small isotopy of $L$ in $S^2$, it is always possible to write $L$ like this). 

Let $\Hc_{\KS}(L)$ denote the stabilized Kauffman-states diagram of $L$ as in the bottom-right corner of Figure~\ref{fig:LumpAndStabilize}, with ladybugs added at each crossing as in Figure~\ref{fig:KSDLWithLadybugs}. Let $\Hc_{\planar}(L)$ denote the planar diagram of $L$. Using the previous constructions, we will give a sequence of Heegaard moves transforming $\Hc_{\KS}(L)$ into $\Hc_{\planar}(L)$. 

\begin{figure}
\includegraphics[scale=0.75]{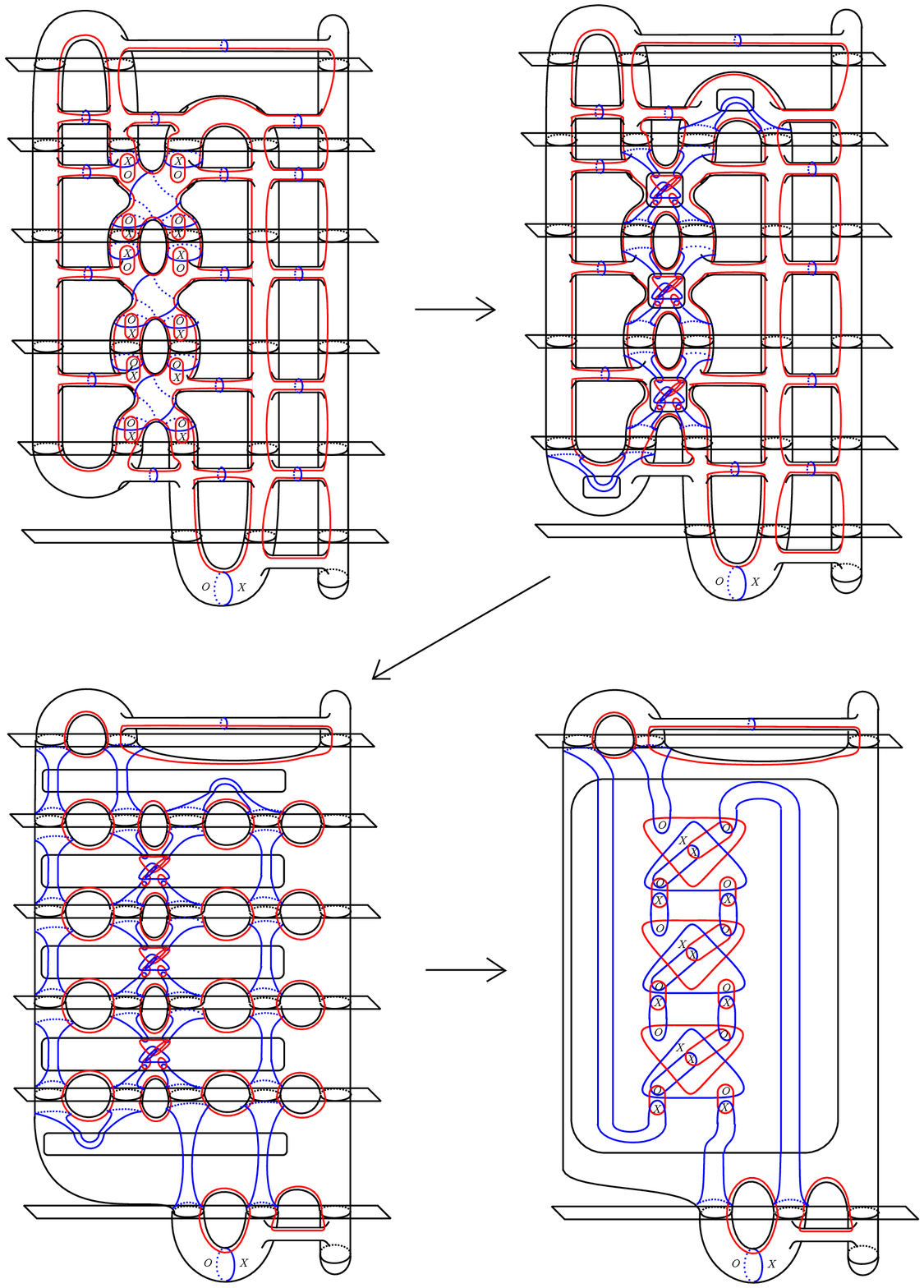}
\caption{Construction~\ref{constr:Main}: part 1. Basepoints are not drawn in the second and third figures, but they should be present as usual.}
\label{fig:MainConst}
\end{figure}

\begin{figure}
\includegraphics[scale=0.75]{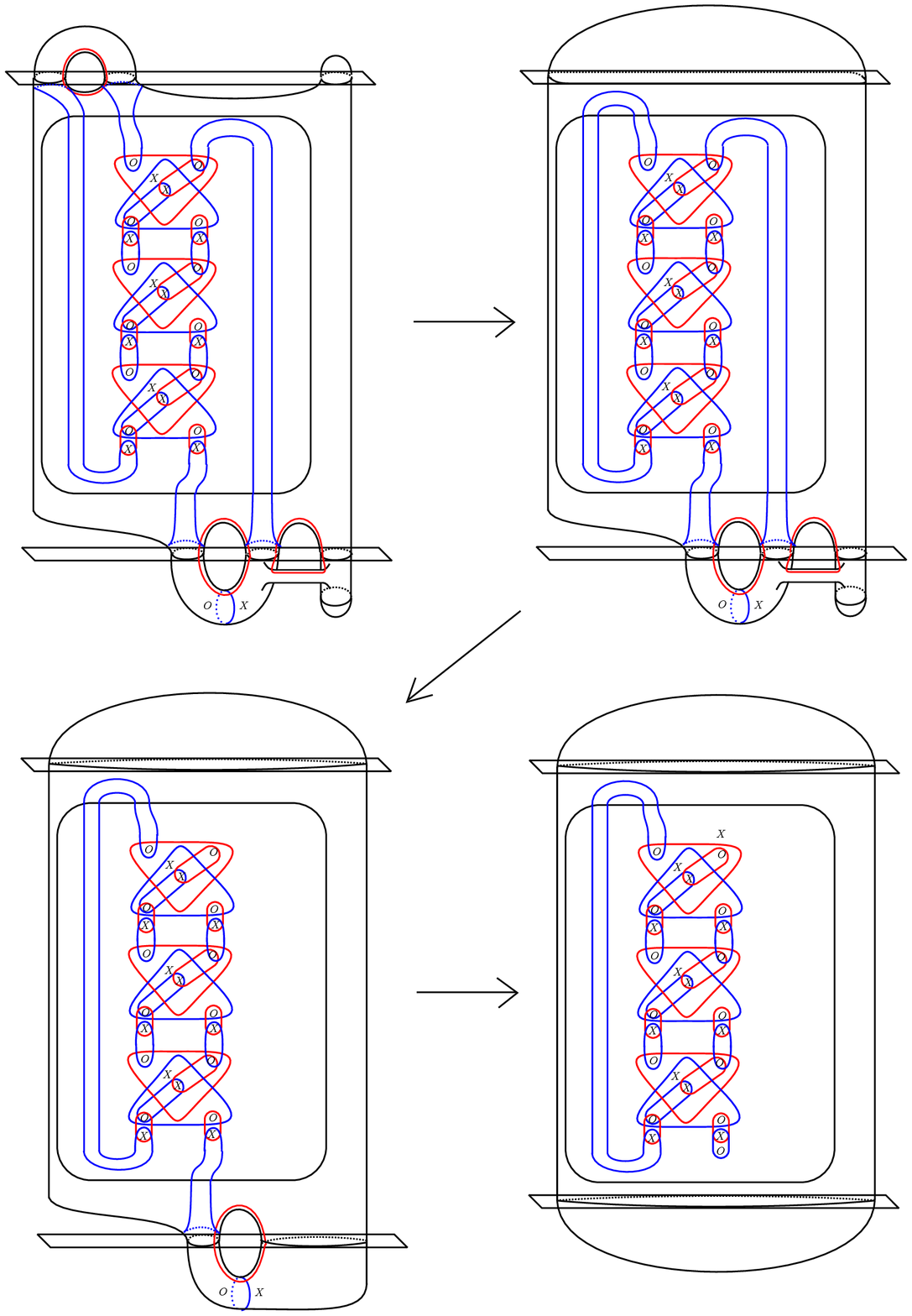}
\caption{Construction~\ref{constr:Main}: part 2.}
\label{fig:MainConstPart2}
\end{figure}

\begin{construction}\label{constr:Main}
Starting from $\Hc_{\KS}(L)$ as shown in the top-left corner of Figure~\ref{fig:MainConst} for the trefoil, use Construction~\ref{constr:Local} at each crossing, maximum point, and minimum point, together with Proposition~\ref{prop:Gluing}, to construct a sequence of Heegaard moves transforming $\Hc_{\KS}(L)$ into a diagram $\Hc'$ like the one shown in the top-right corner of Figure~\ref{fig:MainConst}. Note that each ``row'' of the corresponding Heegaard diagram has the form specified as the starting point of Construction~\ref{constr:HorizExt}. Use Construction~\ref{constr:HorizExt} in each row together with Proposition~\ref{prop:Gluing} to transform $\Hc'$ into a diagram $\Hc''$ like the one shown in the bottom-left corner of Figure~\ref{fig:MainConst}. Then use Construction~\ref{constr:VertExt} and Proposition~\ref{prop:Gluing} to glue the rows together, obtaining a diagram $\Hc'''$ like the one shown in the bottom-right corner of Figure~\ref{fig:MainConst}. 

From this diagram, perform an index one/two destabilization, a handleslide of a $\beta$ circle, and another index one/two destabilization to reach the top-left diagram of Figure~\ref{fig:MainConstPart2}. Perform an index one/two destabilization to reach the bottom-left diagram of Figure~\ref{fig:MainConstPart2}. Finally, slide the basepoints around and perform a final index one/two destabilization to reach the bottom-right diagram of Figure~\ref{fig:MainConstPart2}. This diagram is related to $\Hc_{\planar}(L)$ by a $\beta$ handleslide and an index zero/three destabilization; such modifications are common with the planar diagram, as mentioned in Remark~\ref{rem:SmallerPlanarDiags}.
\end{construction}

\bibliographystyle{plain}
\bibliography{biblio}

\end{document}